\documentclass[11pt,letterpaper]{article}

\usepackage[margin=1in]{geometry}
\usepackage[english]{babel}
\usepackage{amssymb}
\usepackage{mathtools}
\usepackage{graphicx}
\usepackage[round]{natbib}
\usepackage{amsthm}
\usepackage{enumitem}
\usepackage[colorlinks=true,allcolors=blue]{hyperref}
\usepackage{subcaption}
\usepackage{wrapfig}
\newcommand{\R}{\mathbb{R}}
\newcommand{\KL}{\text{KL}}
\DeclareMathOperator*{\argmin}{argmin}
\DeclareMathOperator{\Dkl}{\mathcal{D}_{\KL}}

\newtheorem{theorem}{Theorem}[section]
\newtheorem{proposition}[theorem]{Proposition}
\newtheorem{corollary}[theorem]{Corollary}
\newtheorem{lemma}[theorem]{Lemma}
\newtheorem{remark}[theorem]{Remark}
\newtheorem{assumption}{Assumption}[section]
\theoremstyle{definition}
\newtheorem{definition}[theorem]{Definition}

\newcommand{\papertitle}{Mutual Information Estimation via Triangular Transport Maps:\\Finite Sample Guarantees and Sample Allocation}

\title{\papertitle}
\author{Fengyi Li\thanks{\href{mailto:fengyil@mit.edu}{fengyil@mit.edu}}}
\date{\today}

\begin{document}
\maketitle

\begin{abstract}
We study finite sample guarantees for transport based mutual information estimation. Mutual information can be estimated by first learning the conditional and marginal densities and then evaluating their log-density ratio by Monte Carlo. Given training samples, we use triangular transport maps to learn these densities and then compute the mutual information estimate on an independent batch of evaluation samples. For conditional density estimation using triangular transport maps, we establish non-asymptotic guarantees over realizable H\"older classes by viewing the conditional density as an ordinary density with respect to the product of the law of the conditioning variable and Lebesgue measure. This allows standard density estimation techniques to be applied directly to the conditional density estimation. These results yield finite sample guarantees for the resulting mutual information estimator and characterize how a fixed sample budget should be allocated between training and evaluation in order to minimize the mean squared error. Numerical experiments support the predicted training and evaluation error rates and illustrate the resulting sample allocation tradeoff.
\end{abstract}

\section{Introduction}

Mutual information (MI) measures statistical dependence and quantifies how
much observing one random variable reduces uncertainty about another. It
plays an important role across science and engineering, including material science~\citep{baptista2022bayesian}, remote
sensing~\citep{MI_remote_sense}, clinical trials~\citep{clinical_application},
and optimal experimental design~\citep{lindley,BOED_review,HuanJMacta2024}.
For random vectors $X\in\mathbb R^{d_x}$ and $Y\in\mathbb R^{d_y}$ with
joint law $\mu_{X,Y}$ admitting a Lebesgue density, write $\pi_X,\pi_Y$
for the marginal densities and $\pi_{Y\mid X},\pi_{X\mid Y}$ for the
conditional densities. Then
\begin{align}
I(X;Y)
&=
\mathbb{E}_{\mu_{X,Y}}
\left[
\log
\frac{\pi_{Y \mid X}(Y \mid X)}
{\pi_Y(Y)}
\right]\notag\\
&=
\mathbb{E}_{\mu_{X,Y}}
\left[
\log
\frac{\pi_{X \mid Y}(X \mid Y)}
{\pi_X(X)}
\right].
\label{eq:MI_symmetry}
\end{align}
In Bayesian inference, $X$ is the parameter and $Y$ the observable, so
these two representations correspond to the likelihood-to-marginal and
posterior-to-prior density ratios, respectively. Although the definition is
simple, estimating MI can be difficult when the relevant densities or
expectations are analytically intractable~\citep{HuanJMacta2024,li2024eig}.

A variety of MI estimators have been proposed. Nested Monte Carlo methods
require likelihood evaluations and can converge slowly at finite sample
sizes~\citep{Rainforth2018OnNM,HuanJMacta2024}. When only samples are
available, commonly used alternatives include $k$-nearest-neighbor
estimators~\citep{knn_MI} and variational methods such as
MINE~\citep{MINE}, NWJ~\citep{Nguyen2010_MI}, and related
bounds~\citep{Poole2019OnVB}. Another class of approaches estimates the
densities appearing in the MI density ratio and then evaluates the resulting
plug-in estimator~\citep{Foster2019VariationalBO}. Flexible implementations
include conditional normalizing flows~\citep{Herrmann,dong2024variational}
and Knothe--Rosenblatt (KR) rearrangements
\citep{Koval2024TractableOE,cui2025subspaceacceleratedmeasuretransport,
li2024eig}.

Transport maps are particularly attractive for this purpose because they
represent complicated target distributions through invertible
transformations to simple reference measures
\citep{Marzouk2016,spantini2018jmlr}. With a triangular KR structure,
conditional distributions are naturally encoded in blocks of the map,
making the construction well suited to MI estimation, which requires both
conditional and marginal densities. Moreover, transport maps can be fitted
directly from samples using likelihood-based procedures. Despite this
flexibility, finite-sample guarantees for the resulting MI estimators remain
limited.

For split-sample transport-based MI estimation, \citet{li2024eig} derived
asymptotic guarantees with $N$ training and $M$ independent evaluation
samples. Under in-class parametric assumptions and asymptotic normality,
the bias is $O(N^{-1})$, while the variance and mean-squared error (MSE) are
$O(M^{-1}+N^{-2})$. For a total budget $T=N+M$, their asymptotically
optimal allocation satisfies $M/N\asymp T^{1/3}$. These results, however,
do not provide finite-sample guarantees and rely on asymptotic properties of
the fitted density parameters. Separately, \citet{wang2022minimax}
established non-asymptotic, minimax-optimal density-estimation rates over
H\"older classes, together with convergence rates for triangular transport
maps. These results provide the tools for finite-sample transport analysis,
but do not directly address conditional transport estimation, MI estimation,
or the allocation of samples between training and evaluation.

\subsection{Contribution}
We establish \emph{finite-sample} guarantees for an estimator that learns
$\pi_{Y\mid X}$ and $\pi_Y$ from $N$ shared training samples and evaluates
their fitted log-density ratio using $M$ independent evaluation samples, complementing the asymptotic results of~\citet{li2024eig}. Our analysis separates density-estimation error from Monte Carlo
fluctuation, converting Kullback-Leibler (KL) divergence bounds into high-probability and MSE
guarantees for MI.

For KR maps, we analyze conditional density estimation by treating the
conditional density as an ordinary density with respect to the product of
the conditioning-variable distribution and Lebesgue measure. This yields
non-asymptotic conditional-density rates over realizable H\"older classes
without requiring the conditioning variable itself to admit a density.
Combining the conditional and marginal density-estimation rates gives
explicit high-probability and MSE guarantees for the resulting MI estimator. We then minimize the MSE upper bound
under $T=N+M$ and show that sublinear training suffices for
$O(T^{-1})$ MSE. Numerical experiments examine the training and evaluation
errors and compare the prescribed allocations.

Section~\ref{sec:background} defines the estimator.
Sections~\ref{sec:main_results}--\ref{subsec:allocation} develop the
finite-sample guarantees and sample-allocation results, followed by numerical
experiments in Section~\ref{sec:experiments} and conclusions in
Section~\ref{sec:discussion}. The appendices provide proofs and additional
experimental details.
\section{Problem setup and preliminaries}\label{sec:background}
Let $X\in\mathcal X\subseteq\R^{d_x}$ have an arbitrary Borel probability
law $\mu_X$, and suppose that, conditional on $X=x$,
$Y\in\mathcal Y\subseteq\R^{d_y}$ admits a density
$\pi_{Y\mid X}(\cdot\mid x)$ with respect to Lebesgue measure. The joint
law can therefore be written as
\begin{align}
    \mu_{X,Y}(dx,dy)
    =
    \mu_X(dx)\pi_{Y\mid X}(y\mid x)\,dy.
    \label{eq:cond_disintegration}
\end{align}
The corresponding marginal law is $\mu_Y(dy)=\pi_Y(y)\,dy$, where
$\pi_Y(y)=\int_{\mathcal X}\pi_{Y\mid X}(y\mid x)\,\mu_X(dx)$.
For densities $p,q$ with respect to a common dominating measure $\nu$,
write $\Dkl(p\|q):=\int p\log(p/q)\,d\nu$ for their KL divergence.
Assuming integrability of the log-density ratio, the mutual information reads
\begin{align}
    I(X;Y)
    &=
    \mathbb E_{\mu_{X,Y}}
    \left[
        \log
        \frac{\pi_{Y\mid X}(Y\mid X)}
             {\pi_Y(Y)}
    \right]
    \label{eq:MI_def}\\
    &=
    \int_{\mathcal X}
    \Dkl\!\left(
        \pi_{Y\mid X}(\cdot\mid x)
        \,\middle\|\,
        \pi_Y
    \right)
    \mu_X(dx).
    \label{eq:cond_mi_definition}
\end{align}
Thus, estimating $I(X;Y)$ amounts to estimating the conditional density
$\pi_{Y\mid X}$ and the marginal density $\pi_Y$. Importantly, this
formulation does not require $X$ itself to admit a Lebesgue density.
We abbreviate $I=I(X;Y)$ throughout.
All logarithms are natural. The symbols $\lesssim$ and $\asymp$ denote
one- and two-sided bounds up to positive constants independent of the
sample sizes.

\subsection{Conditional transport representation}
\label{subsec:transport_specialization}

We represent the conditional and marginal densities using triangular
transport maps. For a map $S$ that pushes a target distribution to a known
reference density $\eta$, the density induced by $S$ is
$q_S(z)=\eta(S(z))|\det\nabla S(z)|$. For the exact population transport
$S_0$ associated with a target density $\pi$, this construction satisfies
$q_{S_0}=\pi$. Throughout, we use the term KR
map to refer to a monotone triangular transport map with a prescribed
variable ordering.

For conditional density estimation, we leave $x$ fixed and transport only
the $y$-coordinates. For convenience, we work on the unit cube $Q=[0,1]$,
obtained by rescaling bounded coordinatewise supports when necessary, and
take the reference density $\eta_Y$ to be supported on $Q^{d_y}$.
Write $y_{1:j}=(y_1,\ldots,y_j)$ and $y_{<j}=y_{1:j-1}$, with
$y_{<1}$ the empty vector.

With the variable ordering $(X,Y)$, a full KR map, when it exists, can be
partitioned blockwise as
\begin{align}
S(x,y)
=
\bigl[
S_X(x),\,S_{Y\mid X}(x,y)
\bigr]^\top,
\label{eq:kr_map}
\end{align}
where the lower block has the triangular form
\[
S_{Y\mid X}(x,y)
=
\bigl[
S_1(x,y_1),
\ldots,
S_{d_y}(x,y)
\bigr]^\top.
\]
Our construction works directly with $S_{Y\mid X}$ and therefore does not
require estimation of the leading block $S_X$ or a density model for $X$. Each component $S_j$ depends only on $(x,y_{1:j})$ and is strictly
increasing in $y_j$ for fixed $(x,y_{<j})$. For each $j$ and fixed
$(x,y_{<j})$, the map
$t\mapsto S_j(x,y_{<j},t)$ is a bijection of $Q$ onto itself.
Consequently, for every fixed $x$,
$y\mapsto S_{Y\mid X}(x,y)$ is a bijection of $Q^{d_y}$ onto itself. The conditional density induced by the conditional KR block is
\begin{align}
q_{S_{Y\mid X}}(y\mid x)
&=
\eta_Y\!\left(S_{Y\mid X}(x,y)\right)
\left|\det\nabla_y S_{Y\mid X}(x,y)\right| \notag\\
&=
\eta_Y\!\left(S_{Y\mid X}(x,y)\right)
\prod_{j=1}^{d_y}
\partial_{y_j}S_j(x,y_{1:j}),
\label{eq:cond_KR}
\end{align}
where the second equality follows from the triangular structure of
$\nabla_y S_{Y\mid X}$. Since $S_{Y\mid X}(x,\cdot)$ maps $Q^{d_y}$
bijectively onto itself, a change of variables gives
$\int_{Q^{d_y}}q_{S_{Y\mid X}}(y\mid x)\,dy=1$ for every $x$.

Given $N$ i.i.d.\ training pairs
$\{(X_i,Y_i)\}_{i=1}^N$ from $\mu_{X,Y}$, we fit a maximum likelihood
estimator over an admissible class $\mathcal S_{Y\mid X}$ of conditional
KR maps:
\begin{align}
\widehat S_{Y\mid X}
\in
\argmin_{S_{Y\mid X}\in\mathcal S_{Y\mid X}}
\left\{
-\frac{1}{N}
\sum_{i=1}^N
\log q_{S_{Y\mid X}}(Y_i\mid X_i)
\right\}.
\label{eq:conditional_transport_mle}
\end{align}
We denote the resulting conditional density by
$\hat q_{Y\mid X}=q_{\widehat S_{Y\mid X}}$. The precise regularity and
monotonicity conditions defining $\mathcal S_{Y\mid X}$ are given in
Section~\ref{subsec:conditional_rate}.

The marginal density $\pi_Y$ is estimated analogously by fitting a
$d_y$-dimensional KR map $\widehat S_Y$ to $\{Y_i\}_{i=1}^N$, yielding
$\hat q_Y=q_{\widehat S_Y}$.
Thus, the same training sample is used to construct both density estimators
appearing in the mutual information estimator.
\subsection{Split-sample estimator}\label{subsec:estimator}

We now use the fitted conditional and marginal densities to estimate mutual
information. Let $D_{\mathrm{train}}=\{(X_i,Y_i)\}_{i=1}^N$ denote the
training sample used to construct $\hat q_{Y\mid X}$ and $\hat q_Y$,
and let $D_{\mathrm{eval}}=\{(X_j',Y_j')\}_{j=1}^M$ be an i.i.d.\
evaluation sample from the same joint law $\mu_{X,Y}$, independent of
$D_{\mathrm{train}}$. The split-sample
estimator is
\begin{align}
    \widehat I_{N,M}
    :=
    \frac1M
    \sum_{j=1}^M
    \left[
        \log \hat q_{Y\mid X}(Y_j'\mid X_j')
        -
        \log \hat q_Y(Y_j')
    \right].
    \label{eq:EIG_estimator}
\end{align}
This is the plug-in analogue of~\eqref{eq:MI_def}, with the densities fitted
using $N$ training samples and the resulting log-density ratio evaluated
using $M$ independent samples.
Write $\hat q=(\hat q_{Y\mid X},\hat q_Y)$ for the fitted density pair;
its dependence on $N$ is implicit unless written as $\hat q_N$.
Unconditional expectations and probabilities average over all relevant
training and evaluation randomness; subscripts $\mathrm{train}$ and
$\mathrm{eval}$ specify one source when needed. In particular,
$\operatorname{MSE}(\widehat I_{N,M})
:=\mathbb E[(\widehat I_{N,M}-I)^2]$.
\section{Generic plug-in guarantees}\label{sec:main_results}

In this section, we study the finite-sample error of the mutual information
estimator in a generic setting, without restricting the conditional and
marginal density estimators to any particular form. The goal is to separate
the error arising from density estimation from that arising from Monte Carlo
evaluation. We first identify the density estimation quantities that enter
the plug-in error, and then combine them with concentration of the independent
evaluation sample. The resulting guarantees apply to arbitrary normalized
conditional and marginal density estimators satisfying the assumptions below.
In Section~\ref{subsec:conditional_rate}, we specialize these results to the
conditional and marginal KR estimators introduced above and establish the
required density estimation bounds. For the fitted pair $\hat q$, define
the conditional and marginal KL divergence
\begin{align}
    \mathcal E_C(\hat q)
    &:=
    \mathbb E_{X\sim\mu_X}\!
    \left[
    \Dkl\!\left(
    \pi_{Y\mid X}(\cdot\mid X)
    \,\middle\|\,
    \hat q_{Y\mid X}(\cdot\mid X)
    \right)
    \right],
    \label{eq:conditional_risk}\\
    \mathcal E_Y(\hat q)
    &:=
    \Dkl(\pi_Y\|\hat q_Y).
    \label{eq:marginal_risk}
\end{align}
These quantities measure the errors of the two density estimators entering the
log-density ratio. Since the estimators are learned from the training sample,
both risks remain random through the training data. As shown below, their
difference determines the conditional bias of the MI estimator, while their
individual high-probability bounds control the magnitude of this bias.
We suppress the argument $\hat q$ in $\mathcal E_C$ and $\mathcal E_Y$
when it is unambiguous.

For a fresh pair $(X,Y)\sim\mu_{X,Y}$ independent of training, define the
fitted log-density ratio and its population moments, whenever finite, by
\begin{align}
    \hat\varphi(x,y)
    &:=\log\hat q_{Y\mid X}(y\mid x)-\log\hat q_Y(y),
    \notag\\
    g(\hat q)&:=\mathbb E[\hat\varphi(X,Y)\mid\hat q],
    \notag\\
    V(\hat q)&:=\operatorname{Var}(\hat\varphi(X,Y)\mid\hat q).
    \label{eq:fitted_log_ratio_moments}
\end{align}

\subsection{Assumptions}\label{subsec:assumptions}
The analysis requires two ingredients. First, a tail condition controls
the Monte Carlo fluctuation of the average fitted log-density ratio. Second,
high-probability KL-risk bounds control the error introduced by estimating the
two densities.

\begin{assumption}[Sub-Gaussianity]\label{ass:subgaussian}
There is a deterministic constant $K\geq0$ such that, conditionally on
the fitted pair, for every $t\in\R$,
\begin{align*}
    \mathbb E\!\left[
        e^{t(\hat\varphi(X,Y)
        -g(\hat q))}
        \,\middle|\,\hat q
    \right]
    \leq e^{t^2K^2/2}.
\end{align*}
\end{assumption}

While Assumption~\ref{ass:subgaussian} controls evaluation error for fixed
fitted densities, we also need to control how accurately those densities
approximate the true conditional and marginal densities.

\begin{assumption}[Density estimation error]\label{ass:KL_bounds}
Let $\varepsilon_C(N,\delta)$ and $\varepsilon_Y(N,\delta)$ denote
deterministic upper bounds on the conditional and marginal KL divergence,
respectively, as functions of the training sample size $N$ and confidence
level $\delta$. We assume that, for every $\delta\in(0,1)$,
\begin{align}
    \Pr\!\left(
    \mathcal E_C(\hat q)>\varepsilon_C(N,\delta)
    \right)
    &\leq \delta,
    \notag\\
    \Pr\!\left(
    \mathcal E_Y(\hat q)>\varepsilon_Y(N,\delta)
    \right)
    &\leq \delta.
    \label{eq:KL_assumption}
\end{align}
The risk bounds are non-increasing in both $N$ and $\delta$.
\end{assumption}

The functions $\varepsilon_C$ and $\varepsilon_Y$ are left unspecified in
this generic analysis. They encode the statistical accuracy of the chosen
density estimators. Uniform upper and positive lower bounds on the estimated
densities imply Assumption~\ref{ass:subgaussian} by Hoeffding's lemma
(Lemma~\ref{lem:subgaussian_sufficient}). In
Section~\ref{subsec:conditional_rate}, we verify these conditions for the KR
estimators and derive explicit choices of $\varepsilon_C$ and
$\varepsilon_Y$.

\subsection{Bias identity and finite-sample bound}\label{subsec:main_bound}

We now connect the density-estimation errors above to the MI estimation error.
The key observation is that, conditional on the fitted densities, the bias of
the plug-in estimator is exactly the difference between the marginal and
conditional KL divergences.

\begin{lemma}[Conditional plug-in bias identity]
\label{lem:conditional_bias_identity}
Suppose $\log\frac{\hat q_{Y\mid X}(Y\mid X)}
{\pi_{Y\mid X}(Y\mid X)}
\in L^1(\mu_{X,Y})$ and $\log\frac{\hat q_Y(Y)}{\pi_Y(Y)}
\in L^1(\mu_Y)$. Then
\begin{align}
    \mathbb E[\widehat I_{N,M}-I\mid\hat q]
    =
    \mathcal E_Y(\hat q)-\mathcal E_C(\hat q).
    \label{eq:conditional_bias_identity}
\end{align}
The identity holds for any Borel probability measure $\mu_X$.
\end{lemma}

Thus, the plug-in error has two distinct sources. The first is the fluctuation
of the finite evaluation sample around its conditional mean, and the second
is the error in the fitted densities themselves. In particular,
\begin{align}
    \widehat I_{N,M}-I
    &=
    \underbrace{
    \widehat I_{N,M}-g(\hat q)
    }_{A_M}
    +
    \underbrace{
    \mathcal E_Y(\hat q)-\mathcal E_C(\hat q)
    }_{B(\hat q)},
    \label{eq:decomposition}
\end{align}
with the bias identity proved in
Appendix~\ref{app:proof:lem:conditional_bias_identity}.
In particular, $B(\hat q)=g(\hat q)-I$.
The dependence of $A_M$ on $N$ is suppressed. Here $A_M$ is a Monte Carlo error controlled by the evaluation sample size
$M$, whereas $B(\hat q)$ is a training-dependent plug-in error controlled by
the density-estimation accuracy and hence by $N$. This decomposition is the
basis of both the finite-sample bound below and the sample-allocation analysis
in Section~\ref{subsec:allocation}.

\begin{theorem}[General non-asymptotic bound]\label{thm:main}
Under Assumptions~\ref{ass:subgaussian}--\ref{ass:KL_bounds}, for every
$\delta\in(0,1)$, with probability at least $1-\delta$,
\begin{align}
    |\widehat I_{N,M}-I|
    &\leq
    K\sqrt{\frac{2\log(6/\delta)}{M}}
    \notag\\
    &\quad+
    \varepsilon_C(N,\delta/3)
    +
    \varepsilon_Y(N,\delta/3).
    \label{eq:main_bound}
\end{align}
\end{theorem}

Theorem~\ref{thm:main} makes the two sources of error explicit. The first
term is the Monte Carlo contribution and decays with the evaluation sample
size $M$. The remaining terms are the conditional and marginal
density-estimation errors and decay with the training sample size $N$.
Consequently, once explicit rates for $\varepsilon_C$ and $\varepsilon_Y$
are available, the theorem immediately yields a finite-sample MI error bound
and exposes the tradeoff between training and evaluation samples. The proof details are given in
Appendix~\ref{app:generic}. The same decomposition also yields a corresponding
mean-squared error bound, 
\begin{align}
    \operatorname{MSE}(\widehat I_{N,M})
    \leq
    \frac{K^2}{M}
    +
    \mathbb E[B(\hat q)^2].
    \label{eq:generic_mse}
\end{align}

\section{Conditional density estimation and MI rates via KR maps}
\label{subsec:conditional_rate}
In this section, we first derive a non-asymptotic convergence rate for a
conditional density estimator constructed using KR maps.
We then show that the marginal density inherits the required regularity and
establish the corresponding marginal density-estimation rate. Finally, we
combine these results with the generic guarantees of
Section~\ref{sec:main_results} to obtain finite-sample convergence bounds for
the resulting transport-based mutual information estimator. For the
conditional density analysis, we adapt the proof strategy of
\citet{wang2022minimax} and restrict our attention to H\"older classes on the
unit cube. Throughout, we set $Q=[0,1]$, $\mathcal X=Q^{d_x}$, and
$\mathcal Y=Q^{d_y}$, with $d_x,d_y\geq1$, and let $d=d_x+d_y$. All
regularity assumptions on the conditional density are imposed jointly in
$(x,y)$.

\subsection{Regularity and the conditional estimator}
\label{subsubsec:cond_setting}
We now specify the regularity assumptions and the conditional transport
class used in the convergence analysis.

\begin{definition}[H\"older class]
\label{def:holder}
For $\alpha>0$, let $m_\alpha=\lceil\alpha\rceil-1$ and
$\beta_\alpha=\alpha-m_\alpha\in(0,1]$. On a compact domain $\Omega$,
write $D^\nu$ for the derivative with multi-index $\nu$, $\|\cdot\|_\infty$
for the supremum norm, and $\|\cdot\|$ for the Euclidean norm. Then
the H\"older class
$\mathcal H^\alpha(\Omega)=C^{m_\alpha,\beta_\alpha}(\Omega)$ consists of
functions with continuous derivatives through order $m_\alpha$ and finite
norm
$\|f\|_{\mathcal H^\alpha}:=
\sum_{|\nu|\leq m_\alpha}\|D^\nu f\|_\infty+
\sum_{|\nu|=m_\alpha}\sup_{z\ne z'}
|D^\nu f(z)-D^\nu f(z')|/\|z-z'\|^{\beta_\alpha}$.
Derivatives extend continuously to the boundary. At integer $\alpha$, derivatives of order $\alpha-1$ are taken
to be Lipschitz.
For a vector-valued function $f=(f_1,\ldots,f_k)$, use the componentwise
norm $\|f\|_{\mathcal H^\alpha}
:=(\sum_{j=1}^k\|f_j\|_{\mathcal H^\alpha}^2)^{1/2}$.
\end{definition}

\begin{assumption}[Conditional regularity]
\label{ass:cond_holder}
The conditional density satisfies
$\|\pi_{Y\mid X}\|_{\mathcal H^\alpha(Q^d)}\leq L_0$ and
$0<c\leq\pi_{Y\mid X}(y\mid x)\leq B<\infty$ for some $\alpha>d/2$,
with $\int_{Q^{d_y}}\pi_{Y\mid X}(y\mid x)\,dy=1$ for every
$x\in Q^{d_x}$.
\end{assumption}

\begin{assumption}[Reference density]
\label{ass:reference}
The reference density is
$\eta_Y(u)=\prod_{j=1}^{d_y}\eta_{Y,j}(u_j)$ on $Q^{d_y}$, where each
factor is a probability density satisfying
$\|\eta_{Y,j}\|_{\mathcal H^\alpha(Q)}\leq L_\eta$ and
$0<c_\eta\leq\eta_{Y,j}\leq B_\eta<\infty$.
We write $G_j(t)=\int_0^t\eta_{Y,j}(u)\,du$ for its cumulative distribution
function (CDF).
\end{assumption}

\begin{definition}[Conditional KR map class]
\label{def:cond_class}
For $L\geq1$ and $0<c_{\mathrm{mon}}\leq1$, let
$\mathcal S_{Y\mid X}(\alpha,L,c_{\mathrm{mon}})$ denote the class of
conditional triangular maps introduced in
Section~\ref{subsec:transport_specialization} such that, for every $j$:
\begin{enumerate}[label=(\roman*)]
    \item $S_j(x,y_{<j},\cdot)$ is an increasing bijection of $Q$ onto
    itself for every fixed $(x,y_{<j})$;
    \item $c_{\mathrm{mon}}\leq\partial_{y_j}S_j\leq L$;
    \item $\|S_j\|_{\mathcal H^\alpha(Q^{d_x+j})}\leq L$ and
    $\|\partial_{y_j}S_j\|_{\mathcal H^\alpha(Q^{d_x+j})}\leq L$.
\end{enumerate}
We write
$\mathcal Q_{Y\mid X}
=\{q_{S_{Y\mid X}}:S_{Y\mid X}\in\mathcal S_{Y\mid X}(\alpha,L,c_{\mathrm{mon}})\}$,
where $q_{S_{Y\mid X}}$ is defined in~\eqref{eq:cond_KR}.
\end{definition}

The regularity imposed on the diagonal derivatives
$\partial_{y_j}S_j$, in addition to that of the map components, ensures
that the induced density remains $\alpha$-H\"older. This follows the
transport density-estimation framework of~\citet{wang2022minimax}. We use the maximum conditional-likelihood estimator defined in
\eqref{eq:conditional_transport_mle}, with
$S_{Y\mid X}\in\mathcal S_{Y\mid X}(\alpha,L,c_{\mathrm{mon}})$, and denote the
resulting density by
$\hat q_{Y\mid X}=q_{\widehat S_{Y\mid X}}$.
Existence and measurability of the minimizer are discussed in
Appendix~\ref{app:conditional}.

\subsection{Non-asymptotic conditional and marginal density rates}
\label{subsubsec:cond_mechanism}
In this section, we establish non-asymptotic density-estimation guarantees
for the conditional and marginal KL divergence $\mathcal E_C$ and $\mathcal E_Y$
appearing in Section~\ref{sec:main_results}.

\begin{theorem}[Conditional MLE rate]
\label{thm:cond_rate}
Under Assumptions~\ref{ass:cond_holder} and~\ref{ass:reference}, choose $L$
sufficiently large and $c_{\mathrm{mon}}>0$ sufficiently small so that
$\mathcal S_{Y\mid X}$ contains the true conditional KR map. Let
$a=2\alpha/(2\alpha+d)$. There exist constants $C_1,C_2>0$, depending only
on $\alpha,d_x,d_y,c,B,L_0,L,c_{\mathrm{mon}},L_\eta,c_\eta,B_\eta$ and not
on $\mu_X$, such that, for every $\delta\in(0,1)$,
\begin{align}
    \Pr\!\left(
    \mathcal E_C(\hat q)
    >
    C_1
    \left[
    N^{-a}
    +
    \frac{\log(1/\delta)}{N}
    \right]
    \right)
    &\leq \delta,
    \label{eq:cond_rate_hp}\\
    \mathbb E[\mathcal E_C(\hat q)]
    &\leq C_2N^{-a}.
    \label{eq:cond_rate_exp}
\end{align}
\end{theorem}

The key reduction is to view conditional density estimation as ordinary
density estimation with respect to the dominating measure
$\lambda(dx,dy)=\mu_X(dx)\otimes dy$. Viewed as functions of $(x,y)$,
$\pi_{Y\mid X}(y\mid x)$ and $q_{S_{Y\mid X}}(y\mid x)$ are probability densities with
respect to $\lambda$, so the conditional MLE coincides with an ordinary
density MLE over the corresponding class of $\lambda$-densities. Let
$\mathcal L(S_{Y\mid X})=-\mathbb E_{\mu_{X,Y}}[\log q_{S_{Y\mid X}}(Y\mid X)]$, and let $S_{0,Y\mid X}$
denote the true conditional KR map. Then
\begin{align}
    &\mathcal L(S_{Y\mid X})-\mathcal L(S_{0,Y\mid X})\notag\\
    &\quad=
    \mathbb E_{X\sim\mu_X}
    \left[
    \Dkl
    \left(
    \pi_{Y\mid X}(\cdot\mid X)
    \,\|\,q_{S_{Y\mid X}}(\cdot\mid X)
    \right)
    \right].
    \label{eq:conditional_population_risk}
\end{align}
At $S_{Y\mid X}=\widehat S_{Y\mid X}$, this excess risk equals
$\mathcal E_C(\hat q)$, the conditional KL divergence used in the generic MI
analysis. The remaining argument adapts the density-MLE proof strategy of
\citet{wang2022minimax}. The induced conditional densities form a uniformly
bounded H\"older class with bracketing entropy of order
$\rho^{-d/\alpha}$ at bracket width $\rho>0$. Since $\alpha>d/2$, the corresponding entropy integral
gives a Hellinger critical radius of order
$N^{-\alpha/(2\alpha+d)}$. Uniform upper and positive lower bounds on the
densities then convert squared Hellinger distance into KL divergence, yielding
Theorem~\ref{thm:cond_rate}. The resulting leading KL rate
$N^{-2\alpha/(2\alpha+d)}$ agrees with the minimax-optimal rate established
by \citet{wang2022minimax} for ordinary $d$-dimensional KR density
estimation. Our conditional result adapts their argument by viewing the
conditional density as a density with respect to $\mu_X\otimes dy$. The
marginal rate below follows directly from their ordinary KR result in
dimension $d_y$. The formal reduction, realizability argument,
and entropy bounds are given in
Appendix~\ref{app:proof:thm:cond_rate}.

We next obtain the corresponding rate for the marginal estimator. No
additional regularity assumption on $\pi_Y$ is needed, since its smoothness
and positivity follow directly from those of $\pi_{Y\mid X}$.

\begin{lemma}[Induced regularity of the marginal]
\label{lem:marginal_free}
Under Assumption~\ref{ass:cond_holder},
$\pi_Y\in\mathcal H^\alpha(Q^{d_y})$,
$\|\pi_Y\|_{\mathcal H^\alpha}\leq L_0$, and
$c\leq\pi_Y\leq B$.
\end{lemma}

Lemma~\ref{lem:marginal_free} allows us to apply the ordinary KR
density-estimation result of \citet{wang2022minimax} directly in dimension
$d_y$. Let $\mathcal S_Y(\alpha,L_Y,c_{\mathrm{mon},Y})$ denote the
$d_y$-dimensional analogue of Definition~\ref{def:cond_class}, with
constants chosen so that it contains the true marginal KR map. The marginal
estimator is fitted to the same training observations $\{Y_i\}_{i=1}^N$.

\begin{lemma}[Marginal MLE rate]
\label{lem:marginal_rate_direct}
Let $\widehat S_Y$ be a maximum-likelihood map over $\mathcal S_Y$
based on $\{Y_i\}_{i=1}^N$, let $\hat q_Y=q_{\widehat S_Y}$ be its
induced density, and let $b=2\alpha/(2\alpha+d_y)$. Under
Assumptions~\ref{ass:cond_holder} and~\ref{ass:reference}, there exist
constants $C_{Y,1},C_{Y,2}>0$ such that, for every $\delta\in(0,1)$,
\begin{align}
    \Pr\!\left(
    \mathcal E_Y(\hat q)
    >
    C_{Y,1}
    \left[
    N^{-b}
    +
    \frac{\log(1/\delta)}{N}
    \right]
    \right)
    &\leq \delta,
    \label{eq:marginal_rate_direct_hp}\\
    \mathbb E[\mathcal E_Y(\hat q)]
    &\leq C_{Y,2}N^{-b}.
    \label{eq:marginal_rate_direct_exp}
\end{align}
\end{lemma}

The proof of Lemma~\ref{lem:marginal_free} and the specialization of the
ordinary KR rate are given in Appendix~\ref{app:conditional}.

Finally, the high-probability bounds above yield the second moments needed
for the MSE analysis. Integrating the corresponding KL tails gives
$\mathbb E[\mathcal E_C(\hat q)^2]\leq C_{C,3}N^{-2a}$ and
$\mathbb E[\mathcal E_Y(\hat q)^2]\leq C_{Y,3}N^{-2b}$
for suitable constants $C_{C,3},C_{Y,3}>0$. The tail-integration argument is given in
Lemma~\ref{lem:kl_second_moment} of
Appendix~\ref{app:conditional_moments}.

\subsection{Consequences for mutual information}
\label{subsec:cond_mi}

We now combine the conditional and marginal density-estimation rates with
the generic results of Section~\ref{sec:main_results}.
Theorem~\ref{thm:cond_rate} and Lemma~\ref{lem:marginal_rate_direct}
provide the bounds $\varepsilon_C$ and $\varepsilon_Y$ in
Assumption~\ref{ass:KL_bounds}, while the uniform upper and positive lower
bounds on the fitted transport densities imply
Assumption~\ref{ass:subgaussian}.

\begin{corollary}[Finite-sample MI bounds]
\label{cor:cond_main}
Under Assumptions~\ref{ass:cond_holder} and~\ref{ass:reference}, fit the
conditional and marginal estimators of Theorem~\ref{thm:cond_rate} and
Lemma~\ref{lem:marginal_rate_direct} using the same $N$ training pairs.
Let $[\underline c,\overline B]$ and
$[\underline c_Y,\overline B_Y]$ denote uniform lower and upper bounds for
the fitted conditional and marginal densities, respectively, and define
$K=\frac12\log\!\bigl(
\overline B\,\overline B_Y/
(\underline c\,\underline c_Y)
\bigr)$. Let
$a=2\alpha/(2\alpha+d_x+d_y)$ and
$b=2\alpha/(2\alpha+d_y)$, so that $a<b$.
For the split-sample estimator $\widehat I_{N,M}$ in
\eqref{eq:EIG_estimator}, the following hold:
\begin{enumerate}
    \item
    For every $\delta\in(0,1)$, with probability at least $1-\delta$,
    \begin{align}
        &|\widehat I_{N,M}-I|
        \leq
        K\sqrt{\frac{2\log(6/\delta)}{M}}\notag\\
        +&C_1
        \left[
        N^{-a}
        +
        \frac{\log(3/\delta)}{N}
        \right]
        +C_{Y,1}
        \left[
        N^{-b}
        +
        \frac{\log(3/\delta)}{N}
        \right].
        \label{eq:cond_main_hp}
    \end{align}

    \item
    There exist constants $C',C''>0$ such that
    \begin{align}
        \mathbb E[
        (\widehat I_{N,M}-I)^2
        ]
        \leq
        \frac{C'}{M}
        +
        C''N^{-2a}.
        \label{eq:cond_mse}
    \end{align}
\end{enumerate}
\end{corollary}

The proof is given in Appendix~\ref{app:proof:cor:cond_main}. Since
$a<b$, the conditional density-estimation error decays more slowly than the
marginal error, so the training contribution to the MI MSE is governed by
$N^{-2a}$. Thus, up to constants,
$\mathbb E[(\widehat I_{N,M}-I)^2]
\lesssim M^{-1}+N^{-2a}$.
The two terms respectively capture the Monte Carlo evaluation error,
controlled by $M$, and the density-estimation error, controlled by $N$.
This tradeoff forms the basis of the sample-allocation analysis in the next
section.

\section{Optimal sample allocation}
\label{subsec:allocation}

The MSE bound~\eqref{eq:cond_mse} separates the effects of the training and
evaluation sample sizes as
$\mathbb E[(\widehat I_{N,M}-I)^2]\leq C'/M+C''N^{-2a}$.
Under a fixed total sample budget $T=N+M$, we optimize this upper bound over
the allocation between training and evaluation samples. Here $T$ counts
samples only; the computational cost of fitting the transport maps is not
modeled.

\begin{corollary}[Allocation minimizing the MSE upper bound]
\label{cor:allocation}
Let $a=2\alpha/(2\alpha+d_x+d_y)\in(1/2,1)$ and $C',C''>0$. Under
$N+M=T$, let $(N^\star,M^\star)$ minimize
$C'/M+C''N^{-2a}$. Then, as $T\to\infty$,
\begin{equation}
    N^\star\asymp T^{2/(2a+1)},
    \qquad
    \frac{M^\star}{T}\to1.
    \label{eq:allocation}
\end{equation}
At this allocation, $\mathbb E[(\widehat I_{N^\star,M^\star}-I)^2]\lesssim T^{-1}$. The smaller training order $N\asymp T^{1/(2a)}$, with $M=T-N$, also
attains the $T^{-1}$ MSE rate, but generally does not minimize the upper
bound.
\end{corollary}

The proof is given in Appendix~\ref{app:proof:cor:allocation}. The
distinction between the two training orders follows from minimizing the
bound rather than simply balancing its two terms. In particular, the
first-order condition is
$C'/(T-N^\star)^2
=2aC''(N^\star)^{-(2a+1)}$.
Since $a>1/2$, this implies $N^\star=o(T)$ and hence
$M^\star/T\to1$. More precisely,
$N^\star\sim(2aC''/C')^{1/(2a+1)}T^{2/(2a+1)}$; thus the multiplicative
constant is problem-dependent, and the asymptotic scaling does not provide
a constant-free finite-budget prescription.

For the setting used in the numerical experiments, $d_x=d_y=1$ and
$\alpha=2$, so $a=2/3$. Corollary~\ref{cor:allocation} then gives the
bound-minimizing training order $N^\star\asymp T^{6/7}$, while the smaller
order $N\asymp T^{3/4}$ is already sufficient to attain the $T^{-1}$ MSE
rate. In both cases, $M/T\to1$. More generally, the assumption
$\alpha>(d_x+d_y)/2$ implies $2a>1$, allowing the $T^{-1}$ rate to be
achieved with a sublinear training budget. For comparison, the parametric
bound $O(M^{-1}+N^{-2})$ of \citet{li2024eig} yields the minimizing training
order $T^{2/3}$. As $\alpha\to\infty$ with the dimensions fixed, $a\to1$
and our exponent $2/(2a+1)$ likewise approaches $2/3$. The $O(T^{-1})$ rate matches the canonical Monte Carlo MSE rate and is
the fastest rate permitted by the $M^{-1}$ evaluation term when its
variance is bounded away from zero.

\section{Numerical experiments}
\label{sec:experiments}

We use numerical experiments to examine the density-estimation rates, the
training--evaluation error decomposition, and the sample-allocation behavior
predicted by the theory. All experiments use the same synthetic distribution
with $d_x=d_y=1$ and $X\sim\operatorname{Unif}[0,1]$. Define
$h(t)=(2t-1)|2t-1|$ and $w(x)=(1+h(x))/2$, so that
$\int_0^1h(t)\,dt=0$ and $\int_0^1w(x)\,dx=1/2$. For $x,y\in[0,1]$, let
$\pi_{Y\mid X}(y\mid x)=1+0.6w(x)h(y)$. Since $\int_0^1h(y)\,dy=0$, this is
a normalized conditional density for every $x$. Because $X$ is uniform,
the joint Lebesgue density is
$\pi_{X,Y}(x,y)=\pi_{Y\mid X}(y\mid x)$, while integration over $x$ gives
$\pi_Y(y)=1+0.3h(y)$. On $[0,1]^2$,
$0.4\leq\pi_{Y\mid X}(y\mid x)\leq1.6$, and the density is $C^{1,1}$ but
not $C^2$, corresponding to $\alpha=2$ under
Definition~\ref{def:holder}. The true value of MI is obtained using numerical integration. Sampling and implementation details are given in
Appendix~\ref{app:experiment_benchmark}.

We fit conditional and marginal transport estimators using the same $N$
training pairs. In the present one-dimensional setting, the derivatives of
the KR maps are parameterized by positive, normalized
quadratic B-splines and fitted by maximum likelihood. Further implementation details are given in
Appendix~\ref{app:experiments}. Recall the decomposition
$\widehat I_{N,M}-I=A_M+B(\hat q)$, where $A_M$ is the Monte Carlo
evaluation error and
$B(\hat q)=\mathcal E_Y(\hat q)-\mathcal E_C(\hat q)$ is the
training-dependent plug-in error. With $V(\hat q)$ defined in
\eqref{eq:fitted_log_ratio_moments}, the independence of the training and evaluation samples gives
$\mathbb E[(\widehat I_{N,M}-I)^2]
=\mathbb E[B(\hat q)^2]+\mathbb E[V(\hat q)]/M$.
We examine these two contributions separately before studying their
allocation under a fixed total budget.

\subsection{Training and evaluation error}
\label{subsec:experiment_components}
For this example, $\alpha=2$ and $d_x=d_y=1$, so
$a=2\alpha/(2\alpha+d_x+d_y)=2/3$ and
$b=2\alpha/(2\alpha+d_y)=4/5$. The theory therefore predicts conditional
and marginal KL rates of $N^{-2/3}$ and $N^{-4/5}$, respectively, and a
squared plug-in error of order $N^{-4/3}$.

\begin{figure*}[t]
    \centering

    \begin{subfigure}[t]{0.235\textwidth}
        \centering
        \includegraphics[width=\linewidth]
        {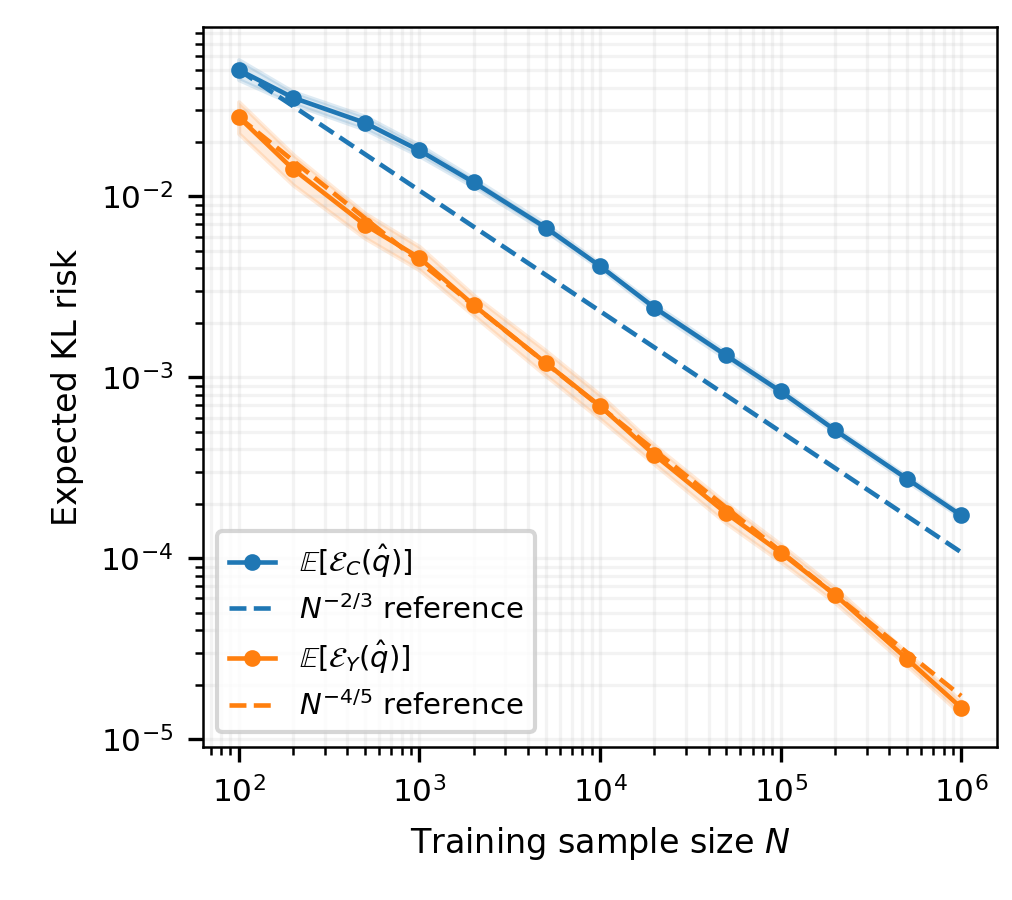}
        \caption{}
        \label{fig:experiment_training_density}
    \end{subfigure}
    \hspace{0.005\textwidth}
    \begin{subfigure}[t]{0.235\textwidth}
        \centering
        \includegraphics[width=\linewidth]
        {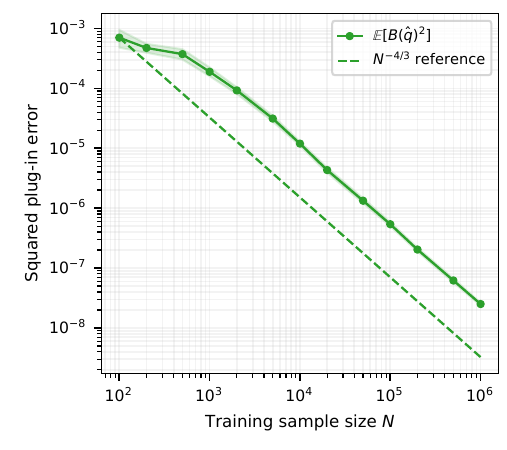}
        \caption{}
        \label{fig:experiment_training_mi}
    \end{subfigure}
    \hspace{0.005\textwidth}
    \begin{subfigure}[t]{0.235\textwidth}
        \centering
        \includegraphics[width=\linewidth]
        {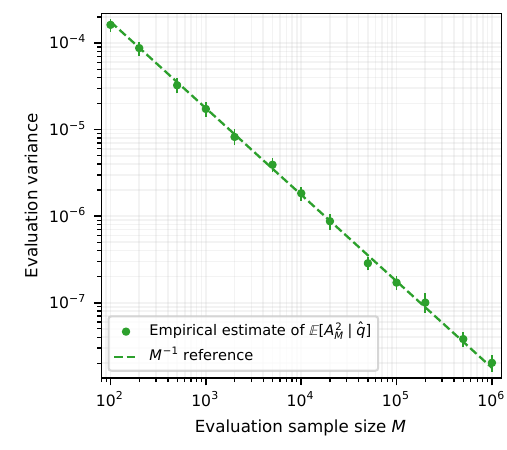}
        \caption{}
        \label{fig:experiment_evaluation_centered}
    \end{subfigure}
    \hspace{0.005\textwidth}
    \begin{subfigure}[t]{0.235\textwidth}
        \centering
        \includegraphics[width=\linewidth]
        {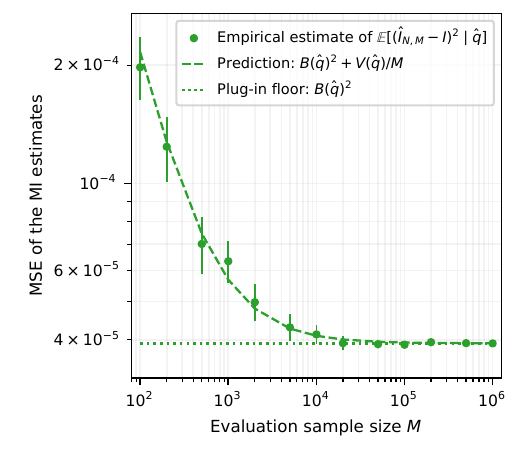}
        \caption{}
        \label{fig:experiment_evaluation_total}
    \end{subfigure}

    \caption{
    Training- and evaluation-sample convergence.
    Panels (a)--(b) vary the training sample size $N$. Panel (a) shows the
    expected conditional and marginal KL divergence, averaged over $30$
    independent training datasets for each $N$, and panel (b) shows the
    corresponding squared plug-in error $\mathbb E[B(\hat q)^2]$.
    Panels (c)--(d) fix the training sample size at $N=5000$, hold the fitted
    densities fixed, and vary the evaluation sample size $M$. Panel (c)
    shows the Monte Carlo evaluation error, while panel (d) shows the MSE
    relative to the true MI and its limiting plug-in floor. Error bars or
    shading denote $95\%$ pointwise bootstrap intervals. 
    }
    \label{fig:experiment_components}
\end{figure*}

For each training sample size $N$, we repeat the fitting procedure on
$30$ independent training datasets. Figures~\ref{fig:experiment_components}(a)
and~\ref{fig:experiment_components}(b) show that the expected conditional
and marginal KL divergence follow the predicted $N^{-2/3}$ and $N^{-4/5}$ rates,
while $\mathbb E[B(\hat q)^2]$ exhibits the predicted $N^{-4/3}$ scaling.
These results support the training-error term $N^{-2a}$ in
\eqref{eq:cond_mse}.

We next fix the training sample size at $N=5000$ and hold the resulting
fitted densities fixed while varying the evaluation sample size $M$.
Conditional on $\hat q$, the evaluation error satisfies
$\mathbb E[A_M^2\mid\hat q]=V(\hat q)/M$, while the MSE
relative to the true MI is
$\mathbb E[(\widehat I_{N,M}-I)^2\mid\hat q]
=B(\hat q)^2+V(\hat q)/M$.
Figures~\ref{fig:experiment_components}(c) and
\ref{fig:experiment_components}(d) confirm the resulting $M^{-1}$ decay
and show that, as $M$ increases, the total error approaches the fixed
plug-in floor $B(\hat q)^2$.

\begin{figure*}[t]
    \centering

    \begin{subfigure}[t]{0.31\textwidth}
        \centering
        \includegraphics[width=\linewidth]
        {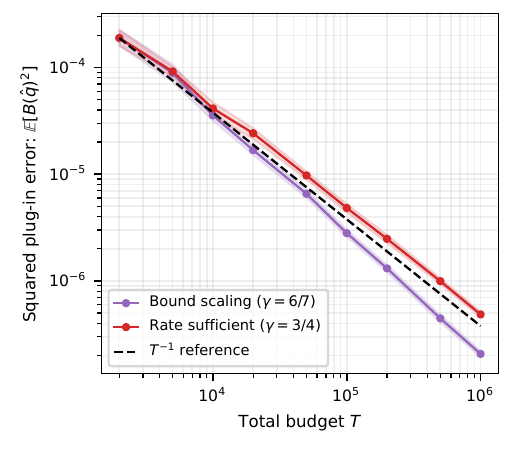}
        \caption{}
        \label{fig:allocation_training}
    \end{subfigure}
    \hspace{0.015\textwidth}
    \begin{subfigure}[t]{0.31\textwidth}
        \centering
        \includegraphics[width=\linewidth]
        {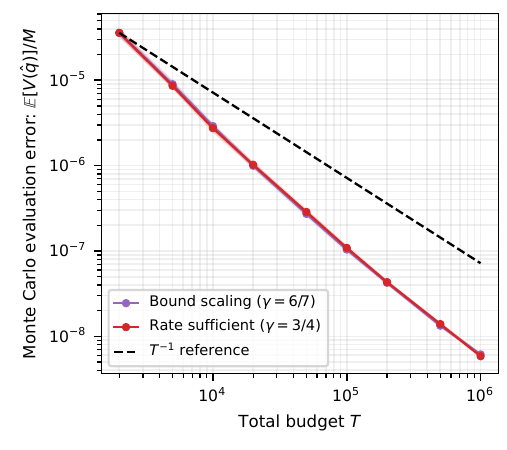}
        \caption{}
        \label{fig:allocation_evaluation}
    \end{subfigure}
    \hspace{0.015\textwidth}
    \begin{subfigure}[t]{0.31\textwidth}
        \centering
        \includegraphics[width=\linewidth]
        {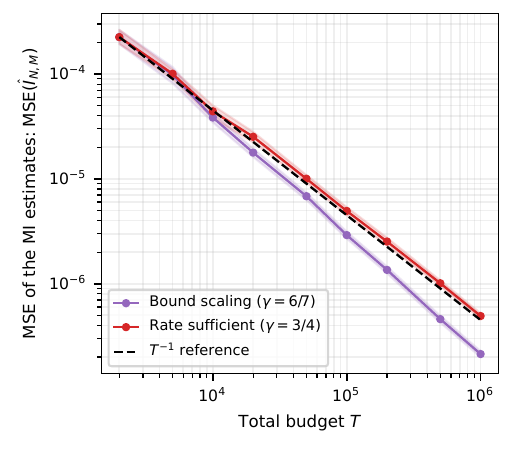}
        \caption{}
        \label{fig:allocation_total}
    \end{subfigure}

    \caption{
    Sample allocation under a fixed total budget. Panel (a) shows the
    empirical squared plug-in error $\mathbb E[B(\hat q)^2]$, panel (b)
    shows the Monte Carlo evaluation contribution
    $\mathbb E[V(\hat q)]/M$, and panel (c) shows the MSE
    $\mathbb E[(\widehat I_{N,M}-I)^2]$ for the
    $N\propto T^{6/7}$ and $N\propto T^{3/4}$ allocations.
    Shading denotes $95\%$ pointwise paired training-bootstrap intervals.
    }
    \label{fig:experiment_allocation}
\end{figure*}

Together, the four panels empirically support the two components of the
finite-sample MSE bound
$\operatorname{MSE}(\widehat I_{N,M})\lesssim N^{-4/3}+M^{-1}$,
which motivates the sample-allocation experiment below.

\subsection{Sample allocation under a fixed budget}
\label{subsec:experiment_allocation}
We next consider a fixed total sample budget $T=N+M$. For each allocation
and each $T$, we repeat the training procedure $R=30$ times using independent
training datasets. For the fitted pair $\hat q^{(r)}$ from repetition $r$,
we compute $B(\hat q^{(r)})$ and $V(\hat q^{(r)})$ by numerical quadrature
against the known joint density. The evaluation-averaged MSE is then
estimated by
\begin{align}
    \widehat R_T(N)
    =
    \frac1R
    \sum_{r=1}^R
    \left[
        B(\hat q^{(r)})^2
        +
        \frac{V(\hat q^{(r)})}{M}
    \right].
    \label{eq:experiment_allocation_risk}
\end{align}
Thus, no evaluation datasets are generated in this experiment; the
expectation over an independent evaluation sample is computed directly by
quadrature. For this example, $a=2/3$. Corollary~\ref{cor:allocation} therefore gives
the upper bound minimizing scaling $N\asymp T^{6/7}$, while
$N\asymp T^{3/4}$ is sufficient to attain the $T^{-1}$ MSE rate. We compare
these two theoretically motivated allocations over nine budgets from
$T=2000$ to $T=10^6$. Specifically, for $\gamma\in\{6/7,3/4\}$, we use
$N_\gamma(T)=\operatorname{round}[1000(T/2000)^\gamma]$,
with $M=T-N_\gamma(T)$. Both schemes are normalized to use \(N=1000\) training samples at \(T=2000\), so they differ only in their scaling exponent. No multiplicative constants are tuned using
the observed risk.

Figure~\ref{fig:experiment_allocation} shows the resulting tradeoff.
The $T^{6/7}$ scaling allocates more samples to training and achieves
a smaller squared plug-in error. Over the tested budget range, its
total MSE also decreases faster than that of the $T^{3/4}$ scaling,
with the two allocations coinciding at $T=2000$. Although the
$T^{3/4}$ scaling leaves more samples for evaluation, its Monte Carlo
contribution is not uniformly smaller because the fitted variance
also depends on the training sample size. The $T^{6/7}$ allocation is obtained by minimizing this finite-sample
upper bound, whereas $T^{3/4}$ is only required to achieve the same
asymptotic $O(T^{-1})$ rate. Accordingly, the theory predicts that the
larger training allocation can yield a smaller finite-budget error even
though both allocations attain the same asymptotic rate. The numerical
results exhibit precisely this tradeoff over the budgets considered.

\section{Conclusion and future work}
\label{sec:discussion}
We developed a finite-sample analysis of transport-based mutual information
estimation using a split-sample procedure. The error is decomposed into a
density-estimation component and a Monte Carlo evaluation component. For KR
maps over realizable H\"older classes, we obtain conditional and marginal KL
rates of
$N^{-2\alpha/(2\alpha+d_x+d_y)}$ and
$N^{-2\alpha/(2\alpha+d_y)}$, respectively, by treating the conditional
density as an ordinary density with respect to the conditioning-variable
distribution and Lebesgue measure. These results yield an MSE bound of order $M^{-1}+N^{-4\alpha/(2\alpha+d_x+d_y)}$ with the conditional density-estimation term dominating the marginal term. Under a fixed budget $T=N+M$, the analysis further characterizes how samples
should be allocated between training and evaluation. The allocation
minimizing the MSE upper bound uses a sublinear number of training samples,
while an even smaller sublinear training allocation is sufficient to attain
the $O(T^{-1})$ MSE rate. Numerical experiments support the predicted training
and evaluation error rates and illustrate the resulting allocation tradeoff.

\section*{AI use statement}
In this work, we used generative AI tools to help develop and refine
theoretical formulations, formulate mathematical claims, assist with the
development and presentation of proofs, provide feedback on research
methodology and numerical experiments. Additionally, we used generative AI tools to assist with writing and editing
parts of the manuscript and improving readability and organization. We take responsibility for the final content of this work, including text,
claims, or artifacts produced with the aid of generative AI.

\bibliographystyle{abbrvnat}
\bibliography{sample}

\clearpage
\appendix
\section{Proofs for the plug-in analysis}\label{app:generic}

\subsection{Proof of Lemma~\ref{lem:conditional_bias_identity}}\label{app:proof:lem:conditional_bias_identity}

\begin{proof}
Conditionally on the fitted densities, independence of the evaluation sample and the disintegration~\eqref{eq:cond_disintegration} give
\begin{align*}
    \mathbb E[\widehat I_{N,M}-I\mid\hat q]
    &=
    \int_{\mathcal X}\int_{\mathcal Y}
    \pi_{Y\mid X}(y\mid x)
    \log\frac{\hat q_{Y\mid X}(y\mid x)}{\pi_{Y\mid X}(y\mid x)}
    \,dy\,\mu_X(dx)\\
    &\quad-
    \int_{\mathcal Y}
    \pi_Y(y)\log\frac{\hat q_Y(y)}{\pi_Y(y)}\,dy,
\end{align*}
which is exactly~\eqref{eq:conditional_bias_identity}. No density of $\mu_X$ with respect to Lebesgue measure is used.
\end{proof}

\subsection{Proof of Theorem~\ref{thm:main}}\label{app:proof:thm:main}

\begin{proof}
Adding and subtracting the conditional mean
$g(\hat q)=\mathbb E[\widehat I_{N,M}\mid\hat q]$, and using
Lemma~\ref{lem:conditional_bias_identity} to identify
$g(\hat q)-I=\mathcal E_Y(\hat q)-\mathcal E_C(\hat q)$, gives
$\widehat I_{N,M}-I=A_M+B(\hat q)$ as in~\eqref{eq:decomposition}.

Conditionally on $D_{\mathrm{train}}$, the fitted log ratio is fixed and the evaluation samples are i.i.d.\ from $\mu_{X,Y}$. For $K>0$, Assumption~\ref{ass:subgaussian} gives
\begin{align*}
    \Pr(|A_M|>t\mid\hat q)\leq2\exp\!\left(-\frac{Mt^2}{2K^2}\right).
\end{align*}
For $K>0$, taking $t=K\sqrt{2\log(6/\delta)/M}$ makes this probability at most $\delta/3$. If $K=0$, the centered summand is zero almost surely and the same conclusion holds directly.

For the plug-in term, Lemma~\ref{lem:conditional_bias_identity} and nonnegativity of KL give
\begin{align}
    |B(\hat q)|\leq\mathcal E_C(\hat q)+\mathcal E_Y(\hat q).
    \label{eq:bias_bound}
\end{align}

To combine these bounds, define
\begin{align*}
    E_1&=\left\{|A_M|\leq K\sqrt{\frac{2\log(6/\delta)}{M}}\right\},\\
    E_2&=\{\mathcal E_C(\hat q)\leq\varepsilon_C(N,\delta/3)\},\\
    E_3&=\{\mathcal E_Y(\hat q)\leq\varepsilon_Y(N,\delta/3)\}.
\end{align*}
Each complement has probability at most $\delta/3$, by the conditional tail bound above and Assumption~\ref{ass:KL_bounds}. Therefore $\Pr(E_1\cap E_2\cap E_3)\geq1-\delta$, without requiring independence between the fitted densities. On this event, the triangle inequality and~\eqref{eq:bias_bound} give~\eqref{eq:main_bound}.
\end{proof}

\subsection{A sufficient condition for sub-Gaussian evaluation}\label{subsec:preliminary}

The following lemma verifies the evaluation-tail assumption from uniform bounds on the fitted densities.

\begin{lemma}[Sufficient condition for sub-Gaussianity]\label{lem:subgaussian_sufficient}
Suppose the estimated densities satisfy the uniform bounds
\begin{align}
    0 < c_\ell \leq \hat{q}_{Y \mid X}(y \mid x) \leq c_u, \qquad 0 < c_\ell' \leq \hat{q}_Y(y) \leq c_u', \label{eq:density_bounds}
\end{align}
for all $(x,y)$ in the support $\mathcal{X} \times \mathcal{Y}$, almost surely over the training data. Then $\hat{\varphi}(X,Y) = \log \hat{q}_{Y \mid X}(Y \mid X) - \log \hat{q}_Y(Y)$ is bounded and Assumption~\ref{ass:subgaussian} holds with
\begin{align*}
    K = \frac12\left[\log(c_u / c_\ell') - \log(c_\ell / c_u')\right]
    = \frac12\log\frac{c_u \, c_u'}{c_\ell \, c_\ell'}.
\end{align*}
\end{lemma}

\begin{proof}
Under~\eqref{eq:density_bounds}, the log-density ratio satisfies $\log(c_\ell / c_u') \leq \hat{\varphi}(x,y) \leq \log(c_u / c_\ell')$ for all $(x,y)$. Centering preserves this interval's length $2K$. Hoeffding's lemma gives the moment-generating-function bound $e^{t^2K^2/2}$ for every $t\in\mathbb R$, as required. For $K>0$ and any $t>0$, it also gives
\begin{align*}
    \Pr\!\left(\left|\frac{1}{M}\sum_{j=1}^{M} \hat{\varphi}(X_j', Y_j') - g(\hat q)\right| > t \;\Bigg|\; \hat{q}\right) \leq 2\exp\!\left(-\frac{Mt^2}{2K^2}\right).
\end{align*}
When $K=0$, the log ratio is constant and the centered average is zero almost surely.
\end{proof}

Note that for the map class of Definition~\ref{def:cond_class}, the positive reference density and uniform monotonicity give a positive lower bound on each induced density. The reference and derivative upper bounds give the corresponding upper bound. Thus Assumption~\ref{ass:subgaussian} follows from the explicitly imposed map-class constraints; smoothness of the true density alone would not suffice.

\section{Proofs for conditional KR estimation}\label{app:conditional}

Throughout this section, we abbreviate the conditional maps
$S_{Y\mid X}$ and $S_{0,Y\mid X}$ by $S$ and $S_0$, respectively,
and write $q_S=q_{S_{Y\mid X}}$.
Thus $S$ denotes only the lower block in this section, not the full map
in~\eqref{eq:kr_map}.

\subsection{Conditional likelihood and population risk}\label{app:conditional_reduction}

The following results justify the dominating-measure reduction used in Section~\ref{subsubsec:cond_mechanism}.

\begin{lemma}[Base-measure reduction]\label{lem:base_measure}
Define the finite measure on $Q^d$
\begin{align}
    \lambda(dx,dy):=\mu_X(dx)\otimes dy.
    \label{eq:base_measure}
\end{align}
Then:
\begin{enumerate}[label=(\roman*)]
    \item for every $S\in\mathcal S_{Y\mid X}$, the function $r_S(x,y):=q_S(y\mid x)$ is a probability density with respect to $\lambda$;
    \item the true conditional density is the $\lambda$-density of $\mu_{X,Y}$, that is,
    \[
        r_0(x,y):=\pi_{Y\mid X}(y\mid x)=\frac{d\mu_{X,Y}}{d\lambda}(x,y);
    \]
    \item the conditional MLE~\eqref{eq:conditional_transport_mle} coincides with the ordinary MLE over the class of $\lambda$-densities $\mathcal R:=\{r_S:S\in\mathcal S_{Y\mid X}\}$ based on the i.i.d.\ sample $(X_i,Y_i)\sim \mu_{X,Y}$.
\end{enumerate}
\end{lemma}

\begin{proof}
Fix $x\in\mathcal X$. By Definition~\ref{def:cond_class}(i)--(ii), the map $y\mapsto S(x,y)$ is a triangular bijection of $Q^{d_y}$ onto itself, strictly increasing in each component's last argument, with Jacobian determinant $\det\nabla_y S=\prod_j\partial_{y_j}S_j>0$. The change of variables $u=S(x,y)$ gives
\begin{align*}
    &\int_{Q^{d_y}}q_S(y\mid x)\,dy\\
    &\quad=\int_{Q^{d_y}}\eta_Y(S(x,y))\det\nabla_y S\,dy\\
    &\quad=\int_{Q^{d_y}}\eta_Y(u)\,du=1.
\end{align*}
Consequently, $r_S\geq0$ and
\begin{align*}
    \int_{Q^d}r_S\,d\lambda
    &=\int_{\mathcal X}\left[\int_{\mathcal Y}q_S(y\mid x)\,dy\right]\mu_X(dx)\\
    &=\int_{\mathcal X}1\,\mu_X(dx)=1.
\end{align*}
Normalization holds for every probability measure $\mu_X$ because the inner integral equals one pointwise in $x$.

For any Borel $A\subseteq Q^d$, the disintegration~\eqref{eq:cond_disintegration} gives
\begin{align*}
    \mu_{X,Y}(A)
    &=\int_{\mathcal X}\int_{\mathcal Y}
    \mathbf 1_A(x,y)\pi_{Y\mid X}(y\mid x)\,dy\,\mu_X(dx)\\
    &=\int_A\pi_{Y\mid X}(y\mid x)\,\lambda(dx,dy).
\end{align*}
Thus $\mu_{X,Y}\ll\lambda$ with density $r_0=\pi_{Y\mid X}$, bounded in $[c,B]$ by Assumption~\ref{ass:cond_holder}. Finally, $\log r_S(X_i,Y_i)=\log q_S(Y_i\mid X_i)$, so the two optimization problems have identical objectives.
\end{proof}

\begin{remark}[The role of $\mu_X$]\label{rem:mu_X_role}
The relevant representation is with respect to $\mu_X\otimes dy$, not Lebesgue measure on $Q^d$. Although $q_S(y\mid x)$ also integrates to one against $dx\,dy$ on the unit cube, that density generally does not describe the observed joint law. With the correct dominating measure, the likelihood and its population risk match the conditional estimation problem. The measure $\mu_X$ is never estimated, and the entropy bound in Lemma~\ref{lem:cond_entropy} is uniform over all probability measures $\mu_X$.
\end{remark}

Under this representation, the excess population risk is exactly the conditional KL divergence.

\begin{lemma}[Exact population identity]\label{lem:cond_population}
Let $\mathcal L(S):=-\mathbb E_{\mu_{X,Y}}[\log q_S(Y\mid X)]$ denote the population risk, and let $S_0\in\mathcal S_{Y\mid X}$ satisfy $q_{S_0}=\pi_{Y\mid X}$. Then, for every $S\in\mathcal S_{Y\mid X}$,
\begin{align}
    &\mathcal L(S)-\mathcal L(S_0)\notag\\
    &\quad=\mathbb E_{X\sim\mu_X}\!\left[
    \Dkl\!\left(\pi_{Y\mid X}(\cdot\mid X)\,\middle\|\,q_S(\cdot\mid X)\right)
    \right]\notag\\
    &\quad=\Dkl(r_0\|r_S)_\lambda,
    \label{eq:cond_population}
\end{align}
where the last expression denotes KL divergence between the $\lambda$-densities $r_0$ and $r_S$.
\end{lemma}

\begin{proof}
Using~\eqref{eq:cond_disintegration},
\begin{align*}
    &\mathcal L(S)-\mathcal L(S_0)\\
    &\quad=\mathbb E_{\mu_{X,Y}}\!\left[\log\frac{\pi_{Y\mid X}(Y\mid X)}{q_S(Y\mid X)}\right]\\
    &\quad=\int_{\mathcal X}\left[
    \int_{\mathcal Y}\pi_{Y\mid X}(y\mid x)
    \log\frac{\pi_{Y\mid X}(y\mid x)}{q_S(y\mid x)}\,dy
    \right]\mu_X(dx).
\end{align*}
The inner integral is the conditional KL divergence. By Lemma~\ref{lem:base_measure}, the same double integral equals $\Dkl(r_0\|r_S)_\lambda$. All quantities are finite because both densities are bounded above and away from zero (Assumption~\ref{ass:cond_holder} and Lemma~\ref{lem:cond_class_regularity}).
\end{proof}

Thus conditional likelihood targets exactly the KL term entering the plug-in bias.

\subsection{Proof of Theorem~\ref{thm:cond_rate}}\label{app:proof:thm:cond_rate}

\paragraph{Existence of the estimator.}
We first establish the existence of a measurable minimizer $\widehat S_{Y\mid X}$. The component functions and their diagonal derivatives form uniformly bounded, equicontinuous families. The H\"older and monotonicity constraints are closed under their joint uniform convergence, and the log-likelihood is continuous because the induced densities are uniformly bounded away from zero. Compactness therefore gives a minimizer, and the measurable maximum theorem supplies a measurable selection. This is the same compactness argument used for the unpenalized H\"older-class MLE in \citet{wang2022minimax}.

Lemmas~\ref{lem:base_measure} and~\ref{lem:cond_population} identify the conditional likelihood with an ordinary density likelihood and its excess population risk with the conditional KL divergence. It remains to establish realizability and entropy bounds, then convert Hellinger concentration into conditional KL control. We begin by verifying that the true conditional density is represented by a map in the prescribed class.

\begin{lemma}[Realizability]\label{lem:cond_realizability}
Under Assumptions~\ref{ass:cond_holder} and~\ref{ass:reference}, there exists $S_0 \in \mathcal S_{Y\mid X}(\alpha, L, c_{\mathrm{mon}})$ with $q_{S_0}(\cdot \mid x) = \pi_{Y \mid X}(\cdot \mid x)$ for every $x$, where
\begin{align*}
    L = L(\alpha, d, c, B, L_0, L_\eta, c_\eta, B_\eta) < \infty, \qquad c_{\mathrm{mon}} = c_{\mathrm{mon}}(c, B, c_\eta, B_\eta) > 0.
\end{align*}
\end{lemma}

\begin{proof}
The construction is the Knothe--Rosenblatt rearrangement applied to $\pi_{Y \mid X}(\cdot \mid x)$ with $x$ carried along as a parameter. Define the partial marginals of the conditional density, for $j = 0, 1, \ldots, d_y$,
\begin{align}
    \bar{p}_j(x, y_{1:j}) := \int_{Q^{d_y - j}} \pi_{Y \mid X}(y_{1:j}, u \mid x)\,du, \label{eq:partial_marginals}
\end{align}
with the conventions $\bar{p}_{d_y} = \pi_{Y \mid X}$ and $\bar{p}_0 \equiv 1$, and set
\begin{align}
    p_j(y_j \mid x, y_{<j}) := \frac{\bar{p}_j(x, y_{1:j})}{\bar{p}_{j-1}(x, y_{1:j-1})}, \qquad j = 1, \ldots, d_y. \label{eq:cond_factors}
\end{align}

\emph{(a) Bounds.}
Since $\pi_{Y \mid X} \in [c, B]$ and the domain of integration in~\eqref{eq:partial_marginals} has unit volume, $c \leq \bar{p}_j \leq B$ for every $j$, and hence
\begin{align}
    0 < \frac{c}{B} \;\leq\; p_j(y_j \mid x, y_{<j}) \;\leq\; \frac{B}{c} < \infty. \label{eq:cond_factor_bounds}
\end{align}

\emph{(b) Regularity.}
Integration over a subset of coordinates does not increase H\"older norms: for any multi-index $\beta$ in the retained variables $(x, y_{1:j})$ with $|\beta| \leq m_\alpha$, differentiation under the integral sign (justified by continuity of the integrand and compactness of the domain) gives $D^\beta \bar{p}_j = \int D^\beta \pi_{Y \mid X}\,du$, whence $\|D^\beta \bar{p}_j\|_\infty \leq \|D^\beta \pi_{Y \mid X}\|_\infty$ and the $\beta_\alpha$-H\"older seminorms of the top-order derivatives obey the same inequality. Thus $\|\bar{p}_j\|_{\mathcal{H}^\alpha} \leq L_0$. Since $\bar{p}_{j-1} \geq c > 0$ and $\mathcal{H}^\alpha(\Omega)$ is an algebra on a bounded convex domain that is stable under division by elements bounded away from zero, \eqref{eq:cond_factors} gives $\|p_j\|_{\mathcal{H}^\alpha} \leq C_a(\alpha, d, c, L_0)$.

\emph{(c) The map.}
Define the conditional KR map by
\begin{align}
    F_j(t; x, y_{<j}) := \int_0^t p_j(u \mid x, y_{<j})\,du, \qquad
    S_{0,j}(x, y_{1:j}) := G_j^{-1}\bigl(F_j(y_j; x, y_{<j})\bigr). \label{eq:cond_KR_construction}
\end{align}
By~\eqref{eq:cond_factor_bounds}, $t \mapsto F_j(t;\cdot)$ is a strictly increasing bijection of $[0,1]$ onto $[0,1]$ with $\partial_t F_j \in [c/B, B/c]$. By Assumption~\ref{ass:reference}, $G_j$ is a strictly increasing bijection of $[0,1]$ onto $[0,1]$ with $G_j' = \eta_{Y,j} \in [c_\eta, B_\eta]$, so $G_j^{-1}$ is a bijection of $[0,1]$ with $(G_j^{-1})' = 1/\eta_{Y,j}(G_j^{-1}) \in [B_\eta^{-1}, c_\eta^{-1}]$ and $\|G_j^{-1}\|_{\mathcal{H}^{\alpha+1}([0,1])} \leq C_b(\alpha, c_\eta, B_\eta, \|\eta_{Y,j}\|_{\mathcal{H}^\alpha})$. It is here that Assumption~\ref{ass:reference} is used in an essential way. Hence $S_{0,j}$ maps $[0,1]$ onto $[0,1]$ in its last argument and is an increasing bijection, giving Definition~\ref{def:cond_class}(i).

\emph{(d) Diagonal derivative.}
By the chain rule and $(G_j^{-1})'(s) = 1/\eta_{Y,j}(G_j^{-1}(s))$,
\begin{align}
    \partial_{y_j} S_{0,j}(x, y_{1:j}) = \frac{p_j(y_j \mid x, y_{<j})}{\eta_{Y,j}\bigl(S_{0,j}(x, y_{1:j})\bigr)} \;\in\; \left[\frac{c}{B B_\eta},\; \frac{B}{c\, c_\eta}\right], \label{eq:cond_diag_deriv}
\end{align}
which gives Definition~\ref{def:cond_class}(ii) with $c_{\mathrm{mon}} = c/(BB_\eta)$. For the regularity in (iii): $F_j \in \mathcal{H}^{\alpha}$ jointly (indeed $\mathcal{H}^{\alpha+1}$ in $y_j$) by the regularity of $p_j$ and integration; composition with the $\mathcal{H}^{\alpha+1}$, Lipschitz function $G_j^{-1}$ preserves $\mathcal{H}^\alpha$; and~\eqref{eq:cond_diag_deriv} exhibits $\partial_{y_j} S_{0,j}$ as a quotient of an $\mathcal{H}^\alpha$ function by $\eta_{Y,j} \circ S_{0,j} \in \mathcal{H}^\alpha$ bounded below by $c_\eta$, hence itself $\mathcal{H}^\alpha$. Collecting the constants yields $S_0 \in \mathcal S_{Y\mid X}(\alpha, L, c_{\mathrm{mon}})$ for $L$ as stated.

\emph{(e) Correctness.}
Finally, by construction and telescoping,
\begin{align*}
    q_{S_0}(y \mid x) =& \prod_{j=1}^{d_y} \eta_{Y,j}\bigl(S_{0,j}\bigr) \cdot \partial_{y_j} S_{0,j}
    = \prod_{j=1}^{d_y} p_j(y_j \mid x, y_{<j})\\
    =& \prod_{j=1}^{d_y} \frac{\bar{p}_j}{\bar{p}_{j-1}}
    = \bar{p}_{d_y}(x,y) = \pi_{Y \mid X}(y \mid x). \qedhere
\end{align*}
\end{proof}

\begin{remark}
Lemma~\ref{lem:cond_realizability} is the conditional analogue of \citet[Proposition~2.1]{wang2022minimax}: positivity and H\"older smoothness of the target density transfer to the KR map through partial integration, division, and composition with the reference CDF. The variable $x$ is a passive parameter throughout, so the construction imposes no regularity on its law.
\end{remark}

To control the complexity of the fitted-density class, we also need uniform bounds for every candidate map.

\begin{lemma}[Uniform regularity of the induced class]\label{lem:cond_class_regularity}
Under Assumption~\ref{ass:reference}, every $q_S \in \mathcal Q_{Y\mid X}$ satisfies
\begin{align}
    0 < \underline{c} := c_\eta^{d_y} c_{\mathrm{mon}}^{d_y} \;\leq\; q_S(y \mid x) \;\leq\; B_\eta^{d_y} L^{d_y} =: \overline{B} < \infty, \label{eq:induced_bounds}
\end{align}
and $\|q_S\|_{\mathcal{H}^\alpha(Q^{d})} \leq C_q$ for a constant $C_q = C_q(\alpha, d, L, c_{\mathrm{mon}}, c_\eta, B_\eta, \|\eta_Y\|_{\mathcal{H}^\alpha})$. Consequently
\begin{align}
    \mathcal Q_{Y\mid X} \;\subseteq\; \mathcal{G} := \bigl\{ g \in \mathcal{H}^\alpha(Q^{d}) : \|g\|_{\mathcal{H}^\alpha} \leq C_q, \;\; \underline{c} \leq g \leq \overline{B} \bigr\}. \label{eq:induced_class}
\end{align}
\end{lemma}

\begin{proof}
The two-sided bound~\eqref{eq:induced_bounds} is immediate from~\eqref{eq:cond_KR}: the product structure in Assumption~\ref{ass:reference} gives $\eta_Y(S(x,y)) \in [c_\eta^{d_y}, B_\eta^{d_y}]$, while each of the $d_y$ diagonal factors lies in $[c_{\mathrm{mon}}, L]$ by Definition~\ref{def:cond_class}(ii). For the H\"older bound: $S \in \mathcal{H}^\alpha$ with values in the compact set $Q^{d_y}$ and is Lipschitz (its entries have bounded first derivatives), so $\eta_Y \circ S \in \mathcal{H}^\alpha$ with norm controlled by $\|\eta_Y\|_{\mathcal{H}^\alpha}$ and $\|S\|_{\mathcal{H}^\alpha} \leq L\sqrt{d_y}$; each factor $\partial_{y_j} S_j \in \mathcal{H}^\alpha$ with norm $\leq L$ by Definition~\ref{def:cond_class}(iii); and $\mathcal{H}^\alpha(Q^d)$ is a Banach algebra on the bounded convex domain $Q^d$, so the product of the $d_y + 1$ factors lies in $\mathcal{H}^\alpha$ with norm bounded by a constant depending only on the listed quantities.
\end{proof}

Lemma~\ref{lem:cond_class_regularity} is the conditional counterpart of \citet[Proposition~2.2]{wang2022minimax}. Its role is to reduce the entropy computation for the \emph{map} class $\mathcal S_{Y\mid X}$ to an entropy computation for a bounded H\"older ball of \emph{densities}, which is classical.

We use this uniform regularity to bound the bracketing entropy. For a class $\mathcal{F}$ of functions and a (semi)norm $\|\cdot\|$, let $N_{[\,]}(\rho, \mathcal{F}, \|\cdot\|)$ denote the $\rho$-bracketing number and $H_{[\,]} = \log N_{[\,]}$ the bracketing entropy.

\begin{lemma}[Bracketing entropy of the induced class]\label{lem:cond_entropy}
Let $\mathcal{G}$ be as in~\eqref{eq:induced_class} and let $\lambda = \mu_X \otimes dy$. Then
\begin{align}
    H_{[\,]}\bigl(\rho, \mathcal{G}, L^2(\lambda)\bigr) \;\leq\; H_{[\,]}\bigl(\rho, \mathcal{G}, \|\cdot\|_\infty\bigr) \;\leq\; A\, \rho^{-d/\alpha}, \qquad 0 < \rho \leq 1, \label{eq:entropy_bound}
\end{align}
where $A = A(\alpha, d, C_q)$ does not depend on $\mu_X$. The same bound, with a different constant, holds for the midpoint square-root class
\begin{align}
    \overline{\mathcal{A}} := \left\{ \sqrt{\frac{g + r_0}{2}} : g \in \mathcal{G} \right\}, \label{eq:sqrt_class}
\end{align}
and there exists a constant $A_1>0$, independent of $\mu_X$, such that
\begin{align}
    H_{[\,]}(\rho,\overline{\mathcal A},L^2(\lambda))
    \leq A_1\rho^{-d/\alpha},
    \qquad 0<\rho\leq1, \label{eq:sqrt_entropy_bound}
\end{align}
\end{lemma}

\begin{proof}
\emph{(a) Sup-norm entropy.}
By the classical entropy estimate of Kolmogorov--Tikhomirov for H\"older balls on a bounded domain \citep[see also][Theorem 2.7.1]{vandervaart1996weak}, the ball $\{g : \|g\|_{\mathcal{H}^\alpha(Q^d)} \leq C_q\}$ satisfies $\log N(\rho, \cdot, \|\cdot\|_\infty) \leq A_0 \rho^{-d/\alpha}$. A $\rho/2$-net in sup norm induces $\rho$-brackets $[g_k - \rho/2,\, g_k + \rho/2]$, so $H_{[\,]}(\rho, \mathcal{G}, \|\cdot\|_\infty) \leq A_0 (\rho/2)^{-d/\alpha}$.

\emph{(b) Transfer to $L^2(\lambda)$, uniformly in $\mu_X$.}
The sup-norm bound transfers to $L^2(\lambda)$ uniformly in $\mu_X$: for any measurable $f, g$ on $Q^d$,
\begin{align*}
    \|f - g\|_{L^2(\lambda)}^2 = \int_{\mathcal{X}} \mu_X(dx) \int_{Q^{d_y}} |f(x,y) - g(x,y)|^2\,dy \;\leq\; \bigl|Q^{d_y}\bigr| \cdot \|f - g\|_\infty^2 \cdot \mu_X(\mathcal{X}) = \|f - g\|_\infty^2,
\end{align*}
using $|Q^{d_y}| = 1$ and $\mu_X(\mathcal{X}) = 1$. Hence any $\rho$-bracket in sup norm is a $\rho$-bracket in $L^2(\lambda)$, and~\eqref{eq:entropy_bound} follows. Crucially, the only property of $\mu_X$ used is that it is a probability measure: no density, smoothness, or positivity enters. This is the precise sense in which the direct conditional argument is free of assumptions on the law of $X$.

\emph{(c) Midpoint square roots.}
By Lemma~\ref{lem:cond_class_regularity} and Assumption~\ref{ass:cond_holder}, every $(g+r_0)/2$ with $g \in \mathcal{G}$ takes values in the compact interval $[(\underline{c}+c)/2,\,(\overline{B}+B)/2] \subset (0,\infty)$ and belongs to a uniformly bounded H\"older ball. Composition with $z \mapsto \sqrt{z}$ therefore places $\overline{\mathcal{A}}$ in a H\"older ball of radius depending only on $(C_q,L_0,\underline{c},\overline{B},c,B,\alpha,d)$. Moreover,
\begin{align*}
    \left\|\sqrt{\frac{g_1+r_0}{2}}-\sqrt{\frac{g_2+r_0}{2}}\right\|_\infty
    \leq \frac{1}{4\sqrt{\underline{c}\wedge c}}\,\|g_1-g_2\|_\infty,
\end{align*}
so the same sup-norm entropy estimate and the preceding norm comparison give~\eqref{eq:sqrt_entropy_bound}.
\end{proof}

The entropy bound determines a critical radius for Hellinger concentration of the MLE.

\begin{lemma}[Entropy integral]\label{lem:entropy_integral}
For $\delta>0$, define the global midpoint square-root class and its entropy integral by
\begin{align}
    \mathcal A
    &:=\left\{\sqrt{\frac{r+r_0}{2}}:r\in\mathcal R\right\},
    \label{eq:global_sqrt_class}\\
    J(\delta)
    &:=\delta+\int_0^\delta
    \sqrt{H_{[\,]}(u,\mathcal A,L^2(\lambda))}\,du.
    \label{eq:global_entropy_integral}
\end{align}
If $\alpha > d/2$, then $J(\delta) < \infty$ for all $0<\delta\leq1$ and
\begin{align}
    J(\delta)
    \leq \delta+\frac{\sqrt{A_1}}{1-d/(2\alpha)}
    \delta^{1-d/(2\alpha)}
    \leq C_{\mathrm{ent}}\delta^{1-d/(2\alpha)},
    \label{eq:entropy_integral}
\end{align}
where $C_{\mathrm{ent}}=1+\sqrt{A_1}/(1-d/(2\alpha))$ is independent of $\mu_X$.
\end{lemma}

\begin{proof}
Since $\mathcal R\subseteq\mathcal G$, we have $\mathcal A\subseteq\overline{\mathcal A}$. Lemma~\ref{lem:cond_entropy} therefore gives
\begin{align*}
    \sqrt{H_{[\,]}(u,\mathcal A,L^2(\lambda))}
    \leq \sqrt{A_1}\,u^{-d/(2\alpha)}.
\end{align*}
The right-hand side is integrable at the origin when $\alpha>d/2$. Direct integration and $\delta\leq\delta^{1-d/(2\alpha)}$ for $0<\delta\leq1$ give~\eqref{eq:entropy_integral}.
\end{proof}

The condition $\alpha>d/2$ is sufficient for this global H\"older-entropy argument. When $\alpha\leq d/2$, the displayed entropy upper bound is not integrable at the origin, so this proof does not yield the stated rate; this observation alone does not imply that the actual entropy integral diverges or that no alternative argument is possible.

\begin{proposition}[MLE concentration at the critical radius]\label{prop:mle_concentration}
Let $\lambda$ be a $\sigma$-finite measure on a standard Borel space $\mathcal Z$, let $\mathcal R$ be a class of $\lambda$-densities containing the true density $r_0>0$ $\lambda$-almost everywhere, and let $\hat r$ be a measurable MLE based on $N$ i.i.d.\ observations. Set $h_\lambda(r,r')=(\tfrac12\int(\sqrt r-\sqrt{r'})^2\,d\lambda)^{1/2}$ and define $J$ from the global class in~\eqref{eq:global_sqrt_class}--\eqref{eq:global_entropy_integral}. If $\delta_N\in(0,1]$ satisfies
\begin{align}
    \sqrt N\,\delta_N^2\geq c_0J(\delta_N),
    \label{eq:critical_radius}
\end{align}
for a suitable absolute constant $c_0$, then there exist absolute constants $c_1,c_2>0$ such that, for every $R\geq\delta_N$,
\begin{align}
    \Pr\bigl(h_\lambda(\hat{r},r_0)>R\bigr)
    \leq c_1\exp(-c_2NR^2). \label{eq:mle_concentration}
\end{align}
\end{proposition}

\begin{proof}
We use a sufficient-condition consequence of \citet[Theorem~7.4]{vandegeer2000empirical}, as reproduced in Theorem~17 of~\citet{kim2016global}. That formulation uses global Hellinger brackets of the midpoint densities, equivalently $L^2$ brackets of their square roots, and the truncated integral
\begin{align*}
    \max\left\{\delta,\,
    \int_{\delta^2/2^{13}}^\delta
    \sqrt{H_{[\,]}(u,\mathcal A,L^2(\lambda))}\,du\right\}
    \leq J(\delta).
\end{align*}
For Lebesgue densities, its conclusion, after converting its Hellinger convention to $h_\lambda$, is bounded by
\begin{align*}
    2^{13/2}\sum_{s=0}^\infty
    \exp\!\left(-\frac{2^{2s}NR^2}{2^{27}}\right),
    \qquad R\geq\delta_N,
\end{align*}
when $c_0=2^{16}$ suffices in~\eqref{eq:critical_radius}. The radius term in $J$ ensures $N\delta_N^2\geq2^{32}$, so this series is at most twice its first term. This gives~\eqref{eq:mle_concentration} with absolute constants.

To justify the extension from Lebesgue densities to the stated dominating measure, let $P_0=r_0\,d\lambda$ and choose a measurable map $\tau:[0,1]\to\mathcal Z$ sending the uniform law to $P_0$; such a map exists because $\mathcal Z$ is standard Borel. For each $r\in\mathcal R$, define the Lebesgue density
\begin{align*}
    q_r(u):=\frac{r(\tau(u))}{r_0(\tau(u))}.
\end{align*}
Then $\int_0^1q_r(u)\,du=1$, $q_{r_0}=1$, and $h_{du}(q_r,1)=h_\lambda(r,r_0)$. For independent uniform $U_i$, the sample $\tau(U_i)$ has law $P_0$, and
\begin{align*}
    \sum_i\log q_r(U_i)
    =\sum_i\log r(\tau(U_i))-\sum_i\log r_0(\tau(U_i)).
\end{align*}
The second term does not depend on $r$, so the fitted density corresponds to an MLE in the transformed class. Brackets also transfer: the order-preserving map $f\mapsto(f/\sqrt{r_0})\circ\tau$ sends the midpoint square roots to $\sqrt{(q_r+1)/2}$ and preserves bracket widths, since
\begin{align*}
    \int_0^1\frac{|f(\tau(u))-g(\tau(u))|^2}{r_0(\tau(u))}\,du
    =\int|f-g|^2\,d\lambda.
\end{align*}
Thus the transformed bracketing entropy is no greater than the original one, and the Lebesgue-density theorem applies. This proves the proposition without assuming a Lebesgue density for $\mu_X$.
\end{proof}

Equivalently, one may bracket $\{\sqrt{(r+r_0)/(2r_0)}:r\in\mathcal R\}$ in $L^2(P_0)$. The division by $\sqrt{r_0}$ is essential for this isometry; raw square-root densities must instead be bracketed in $L^2(\lambda)$. No localization is needed for the global H\"older bound used here.

\begin{lemma}[Solving for $\delta_N$]\label{lem:solve_delta}
Under the conditions of Theorem~\ref{thm:cond_rate}, there are constants $D\geq1$ and $N_0$, independent of $\mu_X$, such that
\begin{align*}
    \delta_N=D N^{-\alpha/(2\alpha+d)}
\end{align*}
belongs to $(0,1]$ and satisfies~\eqref{eq:critical_radius} for every $N\geq N_0$.
\end{lemma}

\begin{proof}
By~\eqref{eq:entropy_integral}, it suffices that $\sqrt N\,\delta_N^{1+d/(2\alpha)}\geq c_0C_{\mathrm{ent}}$. Substituting the stated radius makes the left-hand side equal to $D^{1+d/(2\alpha)}$. Choose $D\geq1$ large enough that this is at least $c_0C_{\mathrm{ent}}$, and then choose $N_0$ so that $\delta_N\leq1$ for $N\geq N_0$. The corresponding power of $N$ is
\begin{align*}
    \frac{1/2}{1+d/(2\alpha)}=\frac{\alpha}{2\alpha+d}.
\end{align*}
\end{proof}

Combining Proposition~\ref{prop:mle_concentration} with Lemma~\ref{lem:solve_delta}, and writing $a=2\alpha/(2\alpha+d)$, yields, for $N\geq N_0$,
\begin{align}
    \Pr\!\left(h_\lambda^2(\hat r,r_0)>t^2D^2N^{-a}\right)
    \leq c_1\exp\!\left(-c_2t^2D^2N^{1-a}\right),
    \qquad t\geq1.
    \label{eq:hellinger_conc}
\end{align}
Here $D$ and $N_0$ depend on the fixed class parameters, while $c_1,c_2$ can be chosen absolute. If the threshold exceeds one, the event is empty.

It remains to convert Hellinger concentration into conditional KL control. Hellinger control does not in general imply KL control---KL can be infinite while Hellinger distance is small. Uniform boundedness of the likelihood ratio, which Lemma~\ref{lem:cond_class_regularity} supplies, closes the gap.

\begin{lemma}[Bounded-ratio KL--Hellinger comparison]\label{lem:kl_hellinger}
Let $p, q$ be probability densities with respect to a measure $\lambda$ such that $m \leq q/p \leq \bar m$ $\lambda$-a.e.\ for constants $0 < m \leq 1 \leq \bar m < \infty$. Then
\begin{align}
    \Dkl(p \,\|\, q) \;\leq\; 2\,C_{m,\bar m}\, h_\lambda^2(p,q), \qquad C_{m,\bar m} := \frac{(1 + \sqrt{\bar m})^2}{2\,m^2}. \label{eq:kl_hellinger}
\end{align}
\end{lemma}

\begin{proof}
Write $t := q/p$ and $\psi(t) := t - 1 - \log t \geq 0$. Because $p$ and $q$ are both probability densities with respect to $\lambda$,
\begin{align}
    \int p\,(t - 1)\,d\lambda = \int (q - p)\,d\lambda = 1 - 1 = 0. \label{eq:centering}
\end{align}
This centering identity is what makes the argument work; it lets us insert the linear term for free:
\begin{align*}
    \Dkl(p \,\|\, q) = \int p \log\frac{p}{q}\,d\lambda = \int p\,(-\log t)\,d\lambda \overset{\eqref{eq:centering}}{=} \int p\,\bigl[t - 1 - \log t\bigr]\,d\lambda = \int p\,\psi(t)\,d\lambda.
\end{align*}
Next we compare $\psi$ with $\phi(t) := (\sqrt{t} - 1)^2$ on $[m, \bar m]$. By Taylor's theorem with Lagrange remainder about $t = 1$, using $\psi(1) = \psi'(1) = 0$ and $\psi''(s) = s^{-2}$, there is $\xi$ between $1$ and $t$ with
\begin{align*}
    \psi(t) = \frac{(t-1)^2}{2\xi^2} \;\leq\; \frac{(t-1)^2}{2 m^2},
\end{align*}
since $\xi \geq \min(1, t) \geq m$. Factoring $(t-1)^2 = (\sqrt{t} - 1)^2 (\sqrt{t} + 1)^2$ and using $t \leq \bar m$,
\begin{align*}
    \psi(t) \;\leq\; \frac{(\sqrt{t}+1)^2}{2m^2}\,\phi(t) \;\leq\; \frac{(1 + \sqrt{\bar m})^2}{2m^2}\,\phi(t) = C_{m,\bar m}\,\phi(t).
\end{align*}
Therefore
\begin{align*}
    \Dkl(p \,\|\, q) &\leq C_{m,\bar m} \int p \left(\sqrt{\tfrac{q}{p}} - 1\right)^{2} d\lambda \\
    &= C_{m,\bar m} \int \bigl(\sqrt{q} - \sqrt{p}\bigr)^2 d\lambda = 2\,C_{m,\bar m}\, h_\lambda^2(p, q). \qedhere
\end{align*}
\end{proof}

\begin{proof}[Proof of Theorem~\ref{thm:cond_rate}]
First take $N\geq N_0$. By Lemma~\ref{lem:base_measure}, $\hat q_{Y\mid X}$ corresponds to the ordinary MLE $\hat r=r_{\widehat S_{Y\mid X}}$ over the class $\mathcal R$ of $\lambda$-densities, and Lemma~\ref{lem:cond_realizability} places the truth $r_0$ in $\mathcal R$. Lemmas~\ref{lem:cond_class_regularity} and~\ref{lem:cond_entropy} bound the global midpoint square-root entropy in $L^2(\lambda)$ uniformly in $\mu_X$. Lemmas~\ref{lem:entropy_integral} and~\ref{lem:solve_delta} and Proposition~\ref{prop:mle_concentration} then give~\eqref{eq:hellinger_conc} with $\delta_N=D N^{-\alpha/(2\alpha+d)}$.

Assumption~\ref{ass:cond_holder} and Lemma~\ref{lem:cond_class_regularity} give the uniform ratio bounds
\begin{align*}
    m := \frac{\underline{c}}{B} \;\leq\; \frac{\hat{r}}{r_0} \;\leq\; \frac{\overline{B}}{c} =: \bar m \qquad \lambda\text{-a.e.},
\end{align*}
almost surely over the training data, so Lemma~\ref{lem:kl_hellinger} applies with these $m, \bar m$ (shrinking $m$ and enlarging $\bar m$ if necessary so that $m \leq 1 \leq \bar m$). Combining with Lemma~\ref{lem:cond_population},
\begin{align}
    \mathbb{E}_{X \sim \mu_X}\!\bigl[\Dkl(\pi_{Y \mid X}(\cdot \mid X) \| \hat{q}_{Y \mid X}(\cdot \mid X))\bigr] = \Dkl(r_0 \| \hat{r})_\lambda \;\leq\; 2 C_{m,\bar m}\, h_\lambda^2(\hat{r}, r_0). \label{eq:cond_kl_to_hellinger}
\end{align}

\emph{(i) High probability.}
Choose
\begin{align*}
    t^2 := \max\left\{1, \; \frac{[\log(c_1/\delta)]_+}{c_2 N \delta_N^2}\right\},
    \qquad [u]_+:=\max\{u,0\}.
\end{align*}
Apply~\eqref{eq:mle_concentration} with $R=t\delta_N$. Its right-hand side is at most $\delta$, so with probability at least $1-\delta$,
\begin{align*}
    h_\lambda^2(\hat r,r_0)
    \leq t^2\delta_N^2
    \leq \delta_N^2+\frac{[\log(c_1/\delta)]_+}{c_2N}.
\end{align*}
Substituting into~\eqref{eq:cond_kl_to_hellinger}, using $\delta_N^2\asymp N^{-2\alpha/(2\alpha+d)}$, and absorbing the constant-order term $N^{-1}[\log c_1]_+$ into $N^{-2\alpha/(2\alpha+d)}$ gives~\eqref{eq:cond_rate_hp}.

\emph{(ii) In expectation.}
Integrating the tail bound~\eqref{eq:hellinger_conc} gives
\begin{align*}
    \mathbb{E}\bigl[h_\lambda^2(\hat{r}, r_0)\bigr]
    &\leq \delta_N^2 + \int_{\delta_N^2}^{\infty} \Pr\bigl(h_\lambda^2 > s\bigr)\,ds
    \;\leq\; \delta_N^2 + \delta_N^2 \int_1^\infty c_1 e^{-c_2 N \delta_N^2 u}\,du \\
    &= \delta_N^2 \left(1 + \frac{c_1 e^{-c_2 N \delta_N^2}}{c_2 N \delta_N^2}\right) \;\leq\; C\,\delta_N^2,
\end{align*}
where the last step uses $N\delta_N^2=D^2N^{d/(2\alpha+d)}\geq1$. Combining with~\eqref{eq:cond_kl_to_hellinger} and taking expectations yields~\eqref{eq:cond_rate_exp}.

Finally, $\mathcal E_C\leq\log(B/\underline c)$ almost surely by the uniform density bounds. Enlarging $C_1,C_2$ to cover the finitely many sample sizes $1\leq N<N_0$ proves both statements for all $N\geq1$. These enlarged constants remain independent of $\mu_X$.
\end{proof}

\subsection{Proof of Lemma~\ref{lem:marginal_free}}\label{app:proof:lem:marginal_free}

\begin{proof}
The regularity argument in the proof of Lemma~\ref{lem:cond_realizability} applies with integration against the probability measure $\mu_X$ in place of Lebesgue measure on a unit-volume set: $D^\beta \pi_Y = \int D^\beta \pi_{Y \mid X}\,\mu_X(dx)$ for $|\beta| \leq m_\alpha$, so all derivative sup-norms and the top-order H\"older seminorm are bounded by those of $\pi_{Y \mid X}$. The pointwise bounds follow since $\mu_X$ is a probability measure.
\end{proof}

\subsection{Proof of Lemma~\ref{lem:marginal_rate_direct}}\label{app:proof:lem:marginal_rate_direct}

\begin{proof}
Lemma~\ref{lem:marginal_free} gives $\pi_Y\in\mathcal H^\alpha(Q^{d_y})$ with the required uniform upper and lower bounds. Apply the realizability, regularity, entropy, and KL--Hellinger arguments of Section~\ref{app:proof:thm:cond_rate} to the ordinary density-estimation problem with dominating measure $dy$ and total dimension $d_y$: realizability and induced-density regularity are the corresponding $d_y$-dimensional KR statements, the H\"older bracketing entropy is $O(\rho^{-d_y/\alpha})$, and Proposition~\ref{prop:mle_concentration} gives the critical radius $N^{-\alpha/(2\alpha+d_y)}$. The KL--Hellinger comparison then yields~\eqref{eq:marginal_rate_direct_hp}--\eqref{eq:marginal_rate_direct_exp}. The entropy integral is finite because $\alpha>(d_x+d_y)/2>d_y/2$.
\end{proof}

\subsection{KL second moments}\label{app:conditional_moments}

\begin{lemma}[Second moments from exponential tails]\label{lem:kl_second_moment}
If a nonnegative random variable $Z$ satisfies $\Pr(Z>u+vt)\leq e^{-t}$ for every $t>0$, with deterministic $u,v>0$, then
\begin{align}
    \mathbb E[Z^2]\leq u^2+2uv+2v^2.
    \label{eq:tail_second_moment}
\end{align}
In particular, the tail bounds in Theorem~\ref{thm:cond_rate} and Lemma~\ref{lem:marginal_rate_direct} imply the second-moment bounds stated in Section~\ref{subsubsec:cond_mechanism}. No independence between the two fitted densities is needed.
\end{lemma}

\begin{proof}
For a nonnegative random variable, Tonelli's theorem gives $\mathbb E[Z^2]=2\int_0^\infty s\,\Pr(Z>s)\,ds$. Splitting at $u$ and substituting $s=u+vt$ in the remaining integral yields
\begin{align*}
    \mathbb E[Z^2]
    &\leq u^2+2v\int_0^\infty(u+vt)e^{-t}\,dt\\
    &=u^2+2uv+2v^2.
\end{align*}
Take $\delta=e^{-t}$ in~\eqref{eq:cond_rate_hp}. The tail condition holds for $Z=\mathcal E_C$, $u=C_1N^{-a}$, and $v=C_1/N$, so
\begin{align*}
    \mathbb E[\mathcal E_C^2]
    &\leq C_1^2\left(N^{-2a}+2N^{-a-1}+2N^{-2}\right)\\
    &\leq5C_1^2N^{-2a},
\end{align*}
using $a<1$ and $N\geq1$. Applying the same argument to~\eqref{eq:marginal_rate_direct_hp}, with $b<1$, proves the marginal bound. No joint event or independence between the two risks is used.
\end{proof}

\subsection{Proof of Corollary~\ref{cor:cond_main}}\label{app:proof:cor:cond_main}

\begin{proof}
The reference bounds and diagonal monotonicity constraints give uniform fitted-density bounds. Lemma~\ref{lem:subgaussian_sufficient} then gives Assumption~\ref{ass:subgaussian} with the stated $K$. Apply Theorem~\ref{thm:main} using the KL tails in Theorem~\ref{thm:cond_rate} and Lemma~\ref{lem:marginal_rate_direct}, each at failure probability $\delta/3$. This proves the high-probability bound.

For the MSE bound, decompose $\widehat I_{N,M}-I=A_M+B(\hat q)$ into the centered evaluation term and the conditional bias. Since $\mathbb E[A_M\mid\hat q]=0$ and $B(\hat q)$ is training-data-measurable, $\mathbb E[A_MB(\hat q)]=0$, so
\begin{align*}
    \mathbb E[(\widehat I_{N,M}-I)^2]
    =\mathbb E[A_M^2]+\mathbb E[B(\hat q)^2].
\end{align*}
Sub-Gaussianity gives $\mathbb E[A_M^2]\leq K^2/M$. Both KL risks are nonnegative, so the bias identity and Lemma~\ref{lem:kl_second_moment} give
\begin{align*}
    \mathbb E[B(\hat q)^2]
    &=\mathbb E[(\mathcal E_Y-\mathcal E_C)^2]\\
    &\leq\mathbb E[\mathcal E_Y^2]+\mathbb E[\mathcal E_C^2]\\
    &\leq C_{Y,3}N^{-2b}+C_{C,3}N^{-2a}\\
    &\leq C''N^{-2a},
\end{align*}
because $b>a$. This proves the MSE bound, without independence between the fitted densities.
\end{proof}

\subsection{Proof of Corollary~\ref{cor:allocation}}\label{app:proof:cor:allocation}

\begin{proof}
Write the displayed upper bound, including its positive constants, as
\begin{align*}
    F_T(M):=\frac{C'}{M}+C''(T-M)^{-2a},
    \qquad 0<M<T.
\end{align*}
It diverges at both endpoints and is strictly convex. Its unique critical point satisfies
\begin{align*}
    \frac{C'}{M^2}=2aC''(T-M)^{-2a-1}.
\end{align*}
Equivalently,
\begin{align*}
    N=\left(\frac{2aC''}{C'}\right)^{1/(2a+1)}M^{2/(2a+1)}.
\end{align*}
Since $a>1/2$, the exponent $2/(2a+1)$ is less than one. Thus $N=o(M)$, $M/T\to1$, and
\begin{align*}
    N\sim\left(\frac{2aC''}{C'}\right)^{1/(2a+1)}T^{2/(2a+1)}.
\end{align*}
At this allocation, $M^{-1}\asymp T^{-1}$ and $N^{-2a}\asymp T^{-4a/(2a+1)}=o(T^{-1})$, proving the MSE claim in Corollary~\ref{cor:allocation}. Alternatively, any allocation satisfying
$N\gtrsim T^{1/(2a)}$ and $T-N\asymp T$ attains the
$O(T^{-1})$ MSE rate. Among polynomially sublinear allocations,
$N\asymp T^{1/(2a)}$ is the smallest sufficient training order.
Since $1/(2a)<2/(2a+1)$ for $a>1/2$, this order is smaller than the
bound-minimizing allocation $N^\star\asymp T^{2/(2a+1)}$.
Rounding to positive integer sample sizes leaves these asymptotic
orders unchanged.
\end{proof}

\section{Experimental details}
\label{app:experiments}

This appendix supplies the implementation details for
the experiments in Section~\ref{sec:experiments}.

\subsection{Fixed finite-smoothness benchmark}
\label{app:experiment_benchmark}

The data-generating conditional, joint, and marginal densities are
specified in Section~\ref{sec:experiments}. In this experimental appendix,
bare $B$, $g$, and $V$ abbreviate $B(\hat q)$, $g(\hat q)$, and $V(\hat q)$
for the fitted pair under discussion.
The generating function $h(t)=(2t-1)|2t-1|$ has continuous first
derivative $h'(t)=4|2t-1|$, which is Lipschitz, while its second
derivative jumps from $-8$ to $8$ at $t=1/2$.
Thus the conditional and marginal densities have the stated
$C^{1,1}$ regularity, with only one curvature breakpoint per coordinate.
Their ranges are respectively $[0.4,1.6]$ and $[0.7,1.3]$.
Writing
\[
 H(y)=\int_0^y h(t)\,dt=\frac{|2y-1|^3-1}{6},
\]
the exact conditional CDF is $y+0.6w(x)H(y)$ and the marginal CDF is
$y+0.3H(y)$.
We sample independent joint pairs by rejection from the uniform
distribution on $[0,1]^2$, with envelope $1.6$. The numerical MI reference is $I=0.0019084$ nats.

\subsection{Spline estimators and numerical controls}
\label{app:experiment_estimators}

Let $\{b_i^{(K)}\}_{i=1}^{K+2}$ denote the quadratic B-spline basis on
$K$ uniform cells with clamped endpoints. The fitted densities are
\begin{align*}
 \hat q_{Y\mid X}(y\mid x)
 &=\sum_{i,j}c_{ij}b_i^{(K_C)}(x)b_j^{(K_C)}(y),&
 K_C&=\lceil2N^{1/6}\rceil,\\
 \hat q_Y(y)&=\sum_jd_jb_j^{(K_Y)}(y),&
 K_Y&=\lceil2N^{1/5}\rceil.
\end{align*}
With $\omega_j^{(K)}=\int_0^1b_j^{(K)}(y)\,dy$, we impose
$\sum_jc_{ij}\omega_j^{(K_C)}=1$ for every row and
$\sum_jd_j\omega_j^{(K_Y)}=1$.
All coefficients lie in $[0.25,1.75]$.
Partition of unity then guarantees the same density bounds.
Linear constraints on derivative coefficients bound every first
partial derivative by 10 and every pure or mixed second partial
derivative by 100, wherever the latter exists. These bounds are fixed
as $N$ increases. Integrating the fitted densities in $y$ gives
increasing conditional and marginal maps from $[0,1]$ onto itself.
The uniform knot grids at different resolutions are not nested, so
approximation error need not change monotonically with the cell count.

We minimize empirical negative log likelihood over these convex
coefficient sets. Conditional fits use SLSQP or constrained Newton steps
with sparse quadratic subproblems solved by OSQP. Marginal fits use SLSQP.
The objective and feasible set are unchanged between the conditional
solvers. A linear-program certificate must give an original-scale
global objective gap at most $10^{-6}$, and actual constraint violation
must not exceed $10^{-7}$.
Failed fits are not silently discarded or replaced, and negative
plug-in MI estimates are not clipped.

Composite Gauss--Legendre quadrature splits at $\{0,1/2,1\}$ and all
fitted knots. Paired KL risks, $g$, and $V$ use orders 24/48, requiring
agreement within $10^{-9}$.
The independent true-MI reference uses orders 48/96.
These are numerical convergence diagnostics, not rigorous integration
error bounds. In Experiment 1, the largest quadrature discrepancy was
$3.47\times10^{-16}$, and the direct residual in
$g-I=\mathcal E_Y-\mathcal E_C$ was at most $2.78\times10^{-17}$.
The finite-dimensional sieve and nonzero optimization tolerance mean
that the experiments are not literal realizations of the
infinite-dimensional MLE in Theorem~\ref{thm:cond_rate}.

\subsection{Sampling and uncertainty}
\label{app:experiment_uncertainty}

\subsubsection{Experiment 1: Training convergence}
\label{subsec:experiment_training}

For each of 13 training sizes from $N=100$ to $N=10^6$, we fit $R=30$
independent density pairs and report
\begin{align}
 \widehat R_C(N)&=\frac1R\sum_{r=1}^R\mathcal E_C(\hat q^{(r)}),
 \notag\\
 \widehat R_Y(N)&=\frac1R\sum_{r=1}^R\mathcal E_Y(\hat q^{(r)}),
 \notag\\
 \widehat R_B(N)&=\frac1R\sum_{r=1}^R[g(\hat q^{(r)})-I]^2.
 \label{eq:experiment_training_risks}
\end{align}
Knot-aware quadrature supplies the population risks and $g$, so no
finite-$M$ evaluation noise enters these measurements. We also check
$g-I=\mathcal E_Y-\mathcal E_C$ for every fitted pair.
The training sizes are $N\in\{100,200,500,1000,2000,5000,10000,20000, 50000,100000,200000,500000,1000000\}$ and the two densities share their training dataset, giving 390 fitted
pairs in total. Matching fits at the same sample size and repetition
are reused across experiments, without selection by observed error.
For each fit, quadrature supplies the vector
$(\mathcal E_C,\mathcal E_Y,B^2)$.
We resample this vector jointly within each $N$, using 2000 bootstrap
replicates, to obtain $95\%$ pointwise percentile intervals for the mean risks.
Log-log slopes are least-squares slopes of log mean risk against log
$N$, using either all 13 sizes or the final decade
$10^5\leq N\leq10^6$. Their bootstrap intervals use the same resampled
mean risks and quantify finite-repetition uncertainty, not uncertainty
in an asymptotic exponent.

\subsubsection{Experiment 2: Fixed-fit evaluation}

For Figures~\ref{fig:experiment_components}(c)--(d), we fix the
$N=5000$ density pair from training repetition 0, chosen without
selection by error. For each $M$ on the same 13-value grid used for $N$
above, we generate 200 independent evaluation datasets, independently
of training and of the datasets at other $M$ values. The densities are
not refitted. With $\widehat I_{N,M}^{(s)}$ denoting evaluation
repetition $s$, the centered empirical risk is
\begin{align}
 \widehat R_A(M\mid\hat q)
 &=\frac1{200}\sum_{s=1}^{200}
   [\widehat I_{N,M}^{(s)}-g(\hat q)]^2.
 \label{eq:experiment_evaluation_risk}
\end{align}
The center $g(\hat q)$ is computed by quadrature, not from these
evaluation repetitions. Replacing $g(\hat q)$ by $I$ gives the empirical
total MSE. Centered and total squared errors from the same repetition
are resampled together, using 2000 bootstrap replicates within each $M$,
to obtain $95\%$ pointwise percentile intervals conditional on this fit.
The reference curves are $V(\hat q)/M$ and
$B(\hat q)^2+V(\hat q)/M$, with the moments computed by quadrature.
The total empirical MSE is measured directly; it need not equal the
empirical centered MSE plus $B(\hat q)^2$ for finitely many repetitions.

\subsubsection{Experiment 3: Prescribed sample allocations}

For Figure~\ref{fig:experiment_allocation}, the total sample sizes are \\
$T\in\{2000,5000,10000,20000,50000,100000,200000,500000,1000000\}$.
We use exactly the two allocation rules specified in
Section~\ref{subsec:experiment_allocation}, with 30 independent training
repetitions at each required $N$. The same fitted pairs are used whenever
the rules request the same $N$; in particular, they coincide at $T=2000$.
Quadrature supplies $(B(\hat q)^2,V(\hat q))$ for each fitted pair.
For each $N$, we jointly resample these pairs 2000 times and compute
the mean squared plug-in error, the mean variance divided by $T-N$,
and their sum as in~\eqref{eq:experiment_allocation_risk}.
The $2.5\%$ and $97.5\%$ bootstrap quantiles give the pointwise bands.
Thus pairing preserves the dependence between the two error components;
no evaluation datasets are generated for this comparison.
Each $T^{-1}$ reference is anchored to the corresponding bound-scaling
component at $T=2000$ and is not a fitted rate or a theorem constant.

\end{document}